\documentclass[11pt,reqno]{siamonline190516}
\usepackage{etoolbox}
\usepackage{geometry}
\usepackage{graphicx}
\usepackage{amssymb}
\usepackage{caption}
\usepackage{subcaption}
\usepackage{enumitem}
\usepackage{mathtools}
\usepackage[all]{xy}
\usepackage{xcolor}
\usepackage{soul}

\usepackage{MnSymbol}

\usepackage[ruled,vlined,algo2e]{algorithm2e}

\SetCommentSty{mycommfont}
\SetKwInput{KwInput}{Input}               
\SetKwInput{KwOutput}{Output}     

\graphicspath{{figures/}}

\numberwithin{equation}{section}

\newsiamremark{remark}{Remark}
\newsiamremark{assumption}{Assumption}

\usepackage[numbers,sort&compress]{natbib}
\bibpunct{[}{]}{;}{n}{,}{,}

\DeclarePairedDelimiterX{\inner}[2]{\langle}{\rangle}{#1, #2}

\DeclareMathOperator*{\ext}{ext}

\allowdisplaybreaks

\title{A Variational Formulation of Accelerated Optimization on Riemannian Manifolds\thanks{Submitted to the editors February 1, 2021.
\funding{This work was supported in part by NSF under grants DMS-1411792, DMS-1345013, DMS-1813635, CCF-2112665, by AFOSR under grant FA9550-18-1-0288, and by the DoD under grant 13106725 (Newton Award for Transformative Ideas during the COVID-19 Pandemic).}}}
\author{Valentin Duruisseaux\thanks{Department of Mathematics, University of California San Diego, La Jolla, CA. (\email{vduruiss@ucsd.edu})} \and Melvin Leok\thanks{Department of Mathematics, University of California San Diego, La Jolla, CA. (\email{mleok@ucsd.edu})}}

\begin{document}

\maketitle

\begin{abstract}
It was shown recently by \cite{SuBoCa16} that Nesterov's accelerated gradient method for minimizing a smooth convex function $f$ can be thought of as the time discretization of a second-order ODE, and that $f(x(t))$ converges to its optimal value at a rate of $\mathcal{O}(1/t^2)$ along any trajectory $x(t)$ of this ODE. A variational formulation was introduced in \cite{WiWiJo16} which allowed for accelerated convergence at a rate of $\mathcal{O}(1/t^p)$, for arbitrary $p>0$, in normed vector spaces. This framework was exploited in~\cite{duruisseaux2020adaptive} using time-adaptive geometric integrators to design efficient explicit algorithms for symplectic accelerated optimization. In \cite{alimisis2020}, a second-order ODE was proposed as the continuous-time limit of a Riemannian accelerated algorithm, and it was shown that the objective function $f(x(t))$ converges to its optimal value at a rate of $\mathcal{O}(1/t^2)$ along solutions of this ODE, thereby generalizing the earlier Euclidean result to the Riemannian manifold setting. In this paper, we show that on Riemannian manifolds, the convergence rate of $f(x(t))$ to its optimal value can also be accelerated to an arbitrary convergence rate $\mathcal{O}(1/t^p)$, by considering a family of time-dependent Bregman Lagrangian and Hamiltonian systems on Riemannian manifolds. This generalizes the results of \cite{WiWiJo16} to Riemannian manifolds and also provides a variational framework for accelerated optimization on Riemannian manifolds. In particular, we will establish results for objective functions on Riemannian manifolds that are geodesically convex, weakly-quasi-convex, and strongly convex. An approach based on the time-invariance property of the family of Bregman Lagrangians and Hamiltonians was used to construct very efficient optimization algorithms in \cite{duruisseaux2020adaptive}, and we establish a similar time-invariance property in the Riemannian setting. This lays the foundation for constructing similarly efficient optimization algorithms on Riemannian manifolds, once the Riemannian analogues of time-adaptive Hamiltonian variational integrators have been developed. The experience with the numerical discretization of variational accelerated optimization flows on vector spaces suggests that the combination of time-adaptivity and symplecticity is important for the efficient, robust, and stable discretization of these variational flows describing accelerated optimization. One expects that a geometric numerical integrator that is time-adaptive, symplectic, and Riemannian manifold preserving will yield a class of similarly promising optimization algorithms on manifolds.
\end{abstract}



\section{Introduction}

Efficient optimization has become one of the major concerns in data analysis. Many machine learning algorithms are designed around the minimization of a loss function or the maximization of a likelihood function. Due to the ever-growing scale of the data sets and size of the problems, there has been a lot of focus on first-order optimization algorithms because of their low cost per iteration. The first gradient descent algorithm was proposed in \cite{Cauchy1847}  by Cauchy to deal with the very large systems of equations he was facing when trying to simulate orbits of celestial bodies, and many gradient-based optimization methods have been proposed since Cauchy's work in 1847. 

In 1983, Nesterov's accelerated gradient method was introduced in \cite{Nes83}, and was shown to converge in $\mathcal{O}(1/k^2)$ to the minimum of the convex objective function $f$, improving on the $\mathcal{O}(1/k)$ convergence rate exhibited by the standard gradient descent methods.
This $\mathcal{O}(1/k^2)$ convergence rate was shown in \cite{Nes04} to be optimal among first-order methods using only information about $\nabla f$ at consecutive iterates. This phenomenon in which an algorithm displays this improved rate of convergence is referred to as acceleration, and other accelerated algorithms have been derived since Nesterov's algorithm, such as accelerated mirror descent \cite{Nem1983} and accelerated cubic-regularized Newton's method \cite{Nes08}. More recently, it was shown in \cite{SuBoCa16} that Nesterov's accelerated gradient method limits to a second-order ODE, as the time-step goes to 0, and that the objective function $f(x(t))$ converges to its optimal value at a rate of $\mathcal{O}(1/t^2)$ along the trajectories of this ODE. It was then shown in \cite{WiWiJo16} that in continuous time, the convergence rate of $f(x(t))$ can be accelerated to an arbitrary convergence rate $\mathcal{O}(1/t^p)$ in normed spaces, by considering flow maps generated by a family of time-dependent Bregman Lagrangian and Hamiltonian systems which is closed under time rescaling. This variational framework and the time-invariance property of the family of Bregman Lagrangians were then exploited in \cite{duruisseaux2020adaptive} using time-adaptive geometric integrators to design efficient explicit algorithms for symplectic accelerated optimization. It was observed that a careful use of adaptivity and symplecticity could result in a significant gain in computational efficiency.  

In the past few years, there has been some effort to derive accelerated optimization algorithms in the Riemannian manifold setting \cite{alimisis2020,Alimisis2020-1,Sra2016, Sra2018, Sra2020, Liu2017}. In \cite{alimisis2020}, a second-order ODE was proposed as the continuous-time limit of a Riemannian accelerated algorithm, and it was shown that the objective function $f(x(t))$ converges to its optimal value at a rate of $\mathcal{O}(1/t^2)$ along solutions of this ODE, generalizing the Euclidean result obtained in \cite{SuBoCa16} to the Riemannian manifold setting. 

In this paper, we show that in continuous time, the convergence rate of $f(x(t))$ to its optimal value can be accelerated to an arbitrary convergence rate $\mathcal{O}(1/t^p)$ on Riemannian manifolds, thereby generalizing the results of \cite{WiWiJo16} to the Riemannian setting. This is achieved by considering a family of time-dependent Bregman Lagrangian and Hamiltonian systems on Riemannian manifolds. This also provides a variational framework for accelerated optimization on Riemannian manifolds, generalizing the normed vector space variational formulation of accelerated optimization introduced in \cite{WiWiJo16}. We will then illustrate the derived theoretical convergence rates by integrating the Bregman Euler--Lagrange equations using a simple numerical scheme to solve eigenvalue and distance minimization problems on Riemannian manifolds. Finally, we will show that the family of Bregman dynamics on Riemannian manifolds is closed under time rescaling, and we will draw inspiration from the approach introduced in \cite{duruisseaux2020adaptive} to take advantage of this invariance property via a carefully chosen Poincar\'e transformation that will allow for the integration of higher-order Bregman dynamics while benefiting from the computational efficiency of integrating lower-order Bregman dynamics on Riemannian manifolds.

\section{Definitions and Preliminaries}

We first introduce the main notions from Riemannian geometry and Lagrangian and Hamiltonian mechanics that will be used throughout this paper (see \cite{MaRa1999,Jost2017,Lee2019,Lang1999,alimisis2020,HaLuWa2006} for more details).

\subsection{Riemannian Geometry}

\begin{definition}
	Given a manifold $\mathcal{Q}$, the \textbf{tangent bundle} $T\mathcal{Q}$ and \textbf{cotangent bundle} $T^*\mathcal{Q}$ are defined by
	\begin{equation*}
	T\mathcal{Q} = \{  (q,v) | q\in \mathcal{Q}, v \in T_q\mathcal{Q}  \}  \qquad \emph{and} \qquad   T^*\mathcal{Q} = \{  (q,p) | q\in \mathcal{Q}, p \in T^*_q\mathcal{Q} \}.
	\end{equation*} 
\end{definition}

\begin{definition} \label{def: fiber metric}
Suppose we have a Riemannian manifold $\mathcal{Q}$ with Riemannian metric  $g(\cdot,\cdot) = \langle \cdot , \cdot \rangle$, represented by the positive-definite symmetric matrix  $(g_{ij}) $ in local coordinates. Then, we define the \textbf{musical isomorphism} $g^{\flat}:T\mathcal{Q} \rightarrow T^*\mathcal{Q}$ by 
	\[ g^{\flat}(u)(v) = g_p(u,v)  \quad  \forall  p\in \mathcal{Q}  \text{ and }  \forall u,v\in T_p\mathcal{Q},  \] and its \textbf{inverse musical isomorphism} $g^{\sharp}:T^*\mathcal{Q} \rightarrow T\mathcal{Q}$. The Riemannian metric $g(\cdot,\cdot) = \langle \cdot , \cdot \rangle$ induces a \textbf{fiber metric} $g^*(\cdot ,\cdot) =   \llangle \cdot , \cdot \rrangle $ on $T^* \mathcal{Q}$ by
	\[  \llangle u , v \rrangle = \langle g^{\sharp}(u), g^{\sharp}(v) \rangle  \quad \forall u,v \in T^* \mathcal{Q}, \]
	represented by the positive definite symmetric matrix $(g^{ij})$ in local coordinates, which is the inverse of the Riemannian metric matrix $(g_{ij}) $. 
\end{definition}

\begin{definition}
	The \textbf{Riemannian gradient} $\emph{gradf}(q) \in T_q  \mathcal{Q}$ at a point $q\in \mathcal{Q}$ of a smooth function $f:\mathcal{Q} \rightarrow \mathbb{R}$ is the tangent vector at $q$ such that 
	\[ \langle  \emph{gradf}(q) , u \rangle = df(q) u   \qquad \forall u\in T_q \mathcal{Q}, \]
	where $df$ is the differential of $f$.
\end{definition}

\begin{definition}
	A \textbf{vector field} on a Riemannian manifold $\mathcal{Q}$ is a map $X:\mathcal{Q} \rightarrow T\mathcal{Q}$ such that $X(q) \in T_q\mathcal{Q}$ for all $q\in \mathcal{Q}$. The set of all vector fields on $\mathcal{Q}$ is denoted $\mathcal{X}(\mathcal{Q})$. The \textbf{integral curve} at $q$ of $X\in \mathcal{X}(\mathcal{Q})$ is the smooth curve $c$ on $\mathcal{Q}$ such that $c(0)=q$ and $c'(t) = X(c(t))$. 
\end{definition}

\begin{definition}
	A \textbf{geodesic} in a Riemannian manifold $\mathcal{Q}$ is a parametrized curve $\gamma : [0,1] \rightarrow \mathcal{Q}$ which is of minimal local length. It can be thought of as a curve having zero ``acceleration" or constant ``speed", that is as a generalization of the notion of straight line from Euclidean spaces to Riemannian manifolds. Given two points $q,\tilde{q} \in \mathcal{Q}$, a vector in $T_q \mathcal{Q}$ can be transported to $T_{\tilde{q}}\mathcal{Q}$ along a geodesic $\gamma$ by an operation $\Gamma(\gamma)_q^{\tilde{q}}:T_q \mathcal{Q}\rightarrow T_{\tilde{q}}\mathcal{Q}$ called \textbf{parallel transport along} $\gamma$. We will simply write $\Gamma_q^{\tilde{q}}$ to denote the parallel transport along some geodesic connecting the two points $q,\tilde{q} \in \mathcal{Q}$, and given $A\in \mathcal{X}(\mathcal{Q})$, we will denote by $\Gamma(A)$ the parallel transport along integral curves of $A$. Note that parallel transport preserves inner products: given a geodesic $\gamma$ from $q\in \mathcal{Q}$ to $\tilde{q}\in \mathcal{Q}$,
	\[ g_q(u,v) = g_{\tilde{q}}\left(\Gamma(\gamma)_{q}^{\tilde{q}} u , \Gamma(\gamma)_{q}^{\tilde{q}} v \right)  \qquad \forall u,v\in T_q \mathcal{Q}.\]
\end{definition}

\begin{definition}
	Given $X,Y \in \mathcal{X}(\mathcal{Q})$, the \textbf{covariant derivative} $\nabla_X Y \in \mathcal{X}(\mathcal{Q})$ of $Y$ along $X$ is
	\[\nabla_X Y(q) = \lim_{h\rightarrow 0} \frac{\Gamma(\gamma)_{\gamma(h)}^{q} Y(\gamma(h)) -Y(q)}{h},\]
	where $\gamma$ is the unique integral curve of $X$ such that $\gamma(0)=q$, for any $q\in \mathcal{Q}$.
\end{definition}

\begin{definition}
	A function $f:\mathcal{Q} \rightarrow \mathbb{R}$ is called $L$\textbf{-smooth} if for any two points $q,\tilde{q} \in \mathcal{Q}$ and geodesic $\gamma$ connecting them,
	\[  \big\| \emph{gradf}(q) - \Gamma(\gamma)_{\tilde{q}}^{q} \emph{gradf}(\tilde{q})  \big\| \leq L \emph{ length}(\gamma).\]
\end{definition}

\begin{definition}
	The \textbf{Riemannian Exponential map} $\emph{Exp}_q:T_q \mathcal{Q} \rightarrow \mathcal{Q}$ at $q\in \mathcal{Q}$ is defined by
	\[ \emph{Exp}_q(v) = \gamma_v(1),  \]
	where $\gamma_v$ is the unique geodesic in $\mathcal{Q}$ such that $\gamma_v(0) = q$ and $\gamma_v'(0) = v$, for any $v\in T_q \mathcal{Q} $. \\  $\emph{Exp}_q$ is a diffeomorphism in some neighborhood $ U \subset T_q\mathcal{Q}$ containing 0, so we can define its inverse map, the \textbf{Riemannian Logarithm map} $\emph{Log}_p : \emph{Exp}_q(U) \rightarrow T_q \mathcal{Q}$.
\end{definition}

\begin{definition}
	Given a Riemannian manifold $\mathcal{Q}$ with sectional curvature bounded below by $K_{\min}$, and an upper bound $D$ for the diameter of the considered domain, define
	\begin{align}\label{eq: zeta}
		\zeta = 
		\begin{cases}
			\sqrt{-K_{\min}} D \coth{ (\sqrt{-K_{\min}} D) }  & \quad \text{if }  K_{\min} < 0 \\ 1   & \quad \text{if }  K_{\min} \geq 0 
		\end{cases} .
	\end{align}
	Note that $\zeta \geq 1$ since $x \coth{x} \geq 1$ for all real values of $x$.
\end{definition}

\subsection{Convexity in Riemannian Manifolds}

\begin{definition}
	A subset $A$ of a Riemannian manifold $\mathcal{Q}$ is called \textbf{geodesically uniquely convex} if every two points of $A$ are connected by a unique geodesic in $A$. A function $f:\mathcal{Q} \rightarrow \mathbb{R}$ is called \textbf{geodesically convex} if for any two points $q,\tilde{q} \in \mathcal{Q}$ and geodesic $\gamma$ connecting them, 
\[ f(\gamma(t)) \leq (1-t) f(q) +t f(\tilde{q})  \qquad \forall t\in [0,1].\] Note that if $f$ is a smooth geodesically convex function on a geodesically uniquely convex subset $A$ of a Riemannian manifold, then
\[ f(q) - f(\tilde{q}) \geq \langle  \emph{gradf}(\tilde{q}) , \emph{Log}_{\tilde{q}}(q) \rangle   \qquad \forall q,\tilde{q} \in A.\]
 A function $f:A\rightarrow \mathbb{R}$ is called \textbf{geodesically $\lambda$-weakly-quasi-convex} with respect to $q \in \mathcal{Q}$ for some $\lambda \in (0,1]$ if
	\[ \lambda \left(f(q) - f(\tilde{q})\right) \geq \langle  \emph{gradf}(\tilde{q}) , \emph{Log}_{\tilde{q}}(q) \rangle   \qquad \forall \tilde{q} \in A.\]
	A function $f:A \rightarrow \mathbb{R}$ is called \textbf{geodesically $\mu$-strongly-convex} for some $\mu >0$ if 
	\[ f(q) - f(\tilde{q}) \geq \langle  \emph{gradf}(\tilde{q}) , \emph{Log}_{\tilde{q}}(q)   \rangle  + \frac{\mu}{2} \| \emph{Log}_{\tilde{q}}{(q)} \|^2  \qquad \forall q,\tilde{q} \in A.\]
	A local minimum of a geodesically convex or weakly-quasi-convex function is also a global minimum, and a geodesically strongly convex function either has no minimum or a unique global minimum. Also note that a geodesically convex function is $\lambda$-weakly-quasi-convex with $\lambda=1$. 
\end{definition}

\subsection{Lagrangian and Hamiltonian Mechanics} 

Given a $n$-dimensional Riemannian manifold $\mathcal{Q}$ with local coordinates $(q^1, \ldots , q^n)$, a \textbf{Lagrangian} is a function $L:T\mathcal{Q} \times \mathbb{R} \rightarrow \mathbb{R}$. The corresponding \textbf{action integral} $\mathcal{S}$ is defined to be the functional
\begin{equation} 
	\mathcal{S} (q) = \int_{0}^{T}{L(q,\dot{q},t)dt},
	\end{equation}
over the space of smooth curves $q:[0,T] \rightarrow \mathcal{Q}$. \textbf{Hamilton's Variational Principle} states that $ \delta S=0$ where the variation $\delta S$ is induced by an infinitesimal variation $\delta q$ of the trajectory $q$ that vanishes at the endpoints. Hamilton's Variational Principle can be shown to be equivalent to the \textbf{Euler--Lagrange equations}
		\begin{equation}\label{eq: EL Basic}
			 \frac{d}{dt} \left( \frac{\partial L} {\partial \dot{q}^k} \right)= \frac{\partial L }{\partial q^k} \qquad \text{for } k=1,\ldots , n.
			\end{equation}
	 The \textbf{Legendre transform} $\mathbb{F}L : T\mathcal{Q} \rightarrow T^* \mathcal{Q}$ of $L$ is defined fiberwise by
	$ \mathbb{F}L : (q^i,\dot{q}^i) \mapsto (q^i,p_i)$ where $p_i = \frac{\partial L}{\partial \dot{q}^i} \in T^* \mathcal{Q}$ is the \textbf{conjugate momentum} of $q^i$. 
	We can then define the associated \textbf{Hamiltonian} $H:T^*\mathcal{Q} \rightarrow \mathbb{R}$ by
\begin{equation}
	H(q,p,t) = \left.\sum_{j=1}^{n} {p_j  \dot{q}^j} - L(q,\dot{q},t)\right|_{p_i=\frac{\partial L}{\partial \dot{q}^i}}.
\end{equation} 
We can also define a Hamiltonian Variational Principle on the Hamiltonian side in momentum phase space
\begin{equation}
	 \delta \int_{0}^{T}{\sum_{j=1}^{n}{\left[ p_j \dot{q}^j - H(q,p,t) \right] dt}} =0 , 
\end{equation}
where the variation is induced by an infinitesimal variation $\delta q$ of the trajectory $q$ that vanishes at the endpoints. This is equivalent to \textbf{Hamilton's equations}, given by
\begin{equation} \dot{p}_k = -\frac{\partial H}{\partial q^k} (p,q),  \qquad  \dot{q}^k = \frac{\partial H}{\partial p_k} (p,q) \qquad \text{for } k=1,\ldots , n,
\end{equation}
which can also be shown to be equivalent to the Euler--Lagrange equations~\eqref{eq: EL Basic}.

\section{Variational Formulation and Convergence Rates}\label{Sec: VariationalFormulation}

\subsection{Inspiration} 

A variational framework was introduced in \cite{WiWiJo16} for accelerated optimization on normed vector spaces. Given a convex, continuously differentiable function $h:\mathcal{X} \rightarrow \mathbb{R}$ on a normed vector space $\mathcal{X}$ such that $\Vert \nabla h(x) \Vert \rightarrow \infty $ as $\Vert x \Vert \rightarrow \infty $, its corresponding Bregman divergence is defined by
\begin{equation}  D_h(x,y) = h(y)-h(x) - \langle \nabla h(x), y-x \rangle .  \end{equation} 
The Bregman Lagrangian and Hamiltonian are then defined to be
\begin{equation} 
	\begin{aligned}
			\mathcal{L}_{\alpha,\beta,\gamma}(x,v,t) & = e^{\alpha_t + \gamma_t} \left[   D_h\left(x+e^{-\alpha_t}v , x\right) - e^{\beta_t} f(x)  \right]  , \\
	\mathcal{H}_{\alpha,\beta,\gamma}(x,r,t) & = e^{\alpha_t + \gamma_t} \left[   D_{h^*}\left(\nabla h(x)+e^{-\gamma_t}r , \nabla h(x) \right) + e^{\beta_t} f(x)  \right] ,
	\end{aligned}
\end{equation}
which are scalar-valued functions of position $x\in \mathcal{X}$, velocity $v\in \mathbb{R}^d$ or  momentum $r\in \mathbb{R}^d$, and of time $t$. Here, $h^*:\mathcal{X}^* \rightarrow \mathbb{R}$ denotes the Legendre transform (or convex dual function) of $h$, defined by $h^*(w) = \sup_{z\in \mathcal{X}}{\left[ \langle w, z \rangle - h(z) \right]}$. The Bregman Lagrangian and Hamiltonian family is parametrized by smooth functions of time, $\alpha_t = \alpha(t),\beta_t = \beta(t),\gamma_t = \gamma(t)$, which are said to satisfy the ideal scaling conditions if
\begin{equation}\label{IdealScaling}
	\dot{\beta}_t \leq e^{\alpha_t} \qquad \text{and} \qquad \dot{\gamma}_t = e^{\alpha_t}.
\end{equation}
If the ideal scaling conditions are satisfied, then by Theorem 1.1 in \cite{WiWiJo16},
\begin{equation}
	f(x(t))-f(x^*) \leq \mathcal{O}(e^{-\beta_t}).
\end{equation}
Another very important property of this family of Bregman Lagrangians is its closure under time dilation, proven in Theorem 1.2 of \cite{WiWiJo16}:
\begin{theorem}\label{ThmTimeDilation}
	If $x(t)$ satisfies the Euler-Lagrange equations corresponding to the Bregman Lagrangian $\mathcal{L}_{\alpha,\beta,\gamma}$, then the reparametrized curve $y(t) = x(\tau(t))$  satisfies the Euler-Lagrange equations corresponding to the modified Bregman Lagrangian $\mathcal{L}_{\tilde{\alpha},\tilde{\beta},\tilde{\gamma}}$ where $\tilde{\alpha}_t = \alpha_{\tau(t)} + \log{\dot{\tau}(t) }$,  $\tilde{\beta}_t = \beta_{\tau(t)}$, and $\tilde{\gamma}_t = \gamma_{\tau(t)}$.
	Furthermore $\alpha,\beta,\gamma$ satisfy the ideal scaling conditions (\ref{IdealScaling}) if and only if $\tilde{\alpha},\tilde{\beta},\tilde{\gamma}$ do.
\end{theorem}

We will now extend these results to the Riemannian manifold setting. Throughout this paper, we will make the following assumptions on the  function $f:\mathcal{Q} \rightarrow \mathbb{R}$ to be minimized and on the ambient Riemannian manifold $\mathcal{Q}$, which are standard assumptions in Riemannian optimization \cite{alimisis2020,Alimisis2020-1,Sra2016,Sra2018}:

\begin{assumption} \label{Assumption 1}
	Solutions of the differential equations derived in this paper remain inside a geodesically uniquely convex subset $A$ of a complete Riemannian manifold $\mathcal{Q}$ (i.e. any two points in $\mathcal{Q}$ can be connected by a geodesic), such that $\emph{diam}(A)$ is bounded above by some constant $D$, that the sectional curvature is bounded from below by $K_{\min}$ on $A$, and that $\emph{Exp}_q$ is well-defined for any $q\in A$, and its inverse $\emph{Log}_q$ is well-defined and differentiable on $A$ for any $q \in A$. Furthermore, $f$ is bounded below, geodesically $L$-smooth and all its minima are inside $A$.  \\
\end{assumption}

\subsection{Convex and Weakly-Quasi-Convex Cases}

Suppose that $f : \mathcal{Q} \rightarrow \mathbb{R}$ is a given geodesically $\lambda$-weakly-quasi-convex function, and that Assumption~\ref{Assumption 1} holds true. Since a geodesically convex function is $\lambda$-weakly-quasi-convex with $\lambda=1$, the following treatment also applies to the case where $f$ is geodesically convex. We define a family of Bregman Lagrangians $\mathcal{L}_{\alpha,\beta,\gamma} :   T\mathcal{Q} \times \mathbb{R} \rightarrow  \mathbb{R} $ parametrized by smooth functions of time $\alpha , \beta, \gamma$ by
\begin{equation} \label{BregmanLGeneral}
	\boxed{\mathcal{L}_{\alpha,\beta,\gamma}(X,V,t) = \frac{1}{2}e^{ \lambda^{-1}\zeta \gamma_t - \alpha_t} \langle V , V\rangle  - e^{\alpha_t + \beta_t +  \lambda^{-1}\zeta \gamma_t} f(X) ,}
\end{equation}
and the corresponding Bregman Hamiltonians $\mathcal{H}_{\alpha,\beta,\gamma} :   T^*\mathcal{Q} \times \mathbb{R} \rightarrow  \mathbb{R} $ are given by 
\begin{equation} \label{BregmanHGeneral}
	\boxed{	\mathcal{H}_{\alpha,\beta,\gamma}(X,R,t)= \frac{1}{2}e^{ \alpha_t - \lambda^{-1} \zeta \gamma_t } \llangle R, R\rrangle  +  e^{\alpha_t + \beta_t + \lambda^{-1} \zeta \gamma_t} f(X) ,}
\end{equation}
\normalsize where $X \in \mathcal{Q}$ denotes position on the manifold $\mathcal{Q}$, $V$ is the velocity vector field, $R$ is the momentum covector field, $t$ is the time variable, and $\zeta$ is given by equation~\eqref{eq: zeta}. This family of functions is a generalization of the Bregman Lagrangians and Hamiltonians introduced in \cite{WiWiJo16} for the convex continuously differentiable function $h(x) = \frac{1}{2}\langle x,x\rangle$. Throughout this paper, we will assume that the parameter functions $\alpha, \beta, \gamma$ satisfy the ideal scaling conditions~\eqref{IdealScaling}.

\begin{theorem} 
	The Bregman Euler--Lagrange equation corresponding to the Bregman Lagrangian $\mathcal{L}_{\alpha,\beta,\gamma}$ is given by
	\begin{equation} \label{eq: EL}
		\boxed{	\nabla_{\dot{X}}\dot{X}  +\left( \lambda^{-1} \zeta e^{\alpha_t} -  \dot{\alpha}_t \right) \dot{X} + e^{2\alpha_t+\beta_t }  \emph{gradf}(X) = 0.}
	\end{equation}
 \proof{See Appendix~\ref{Section: Appendix Derivation EL}.}
\end{theorem} 
\begin{theorem}\label{Thm: ConvergenceRateC}
	Suppose that $f : \mathcal{Q} \rightarrow \mathbb{R}$ is a geodesically $\lambda$-weakly-quasi-convex function, and that Assumption~\ref{Assumption 1} is satisfied. Then, any solution $X(t)$ to the Bregman Euler--Lagrange equation~\eqref{eq: EL} converges to a minimizer $x^*$ of $f$ with rate	\begin{equation}
		\boxed{f(X(t)) - f(x^*) \leq \frac{2\lambda^2 e^{\beta_0} \left( f(x_0) - f(x^*) \right) +  \zeta \| \emph{Log}_{x_0}{(x^*)} \|^2 }{2\lambda^2 e^{\beta_t}} = \mathcal{O}(e^{-\beta_t}).}
	\end{equation} 
\proof{ See Appendix~\ref{Section: Appendix Proof Convergence Rate}. }
\end{theorem}    

A $p>0$ parametrized subfamily of Bregman Lagrangians and Hamiltonians, that is of particular practical interest, is given by the choice of parameter functions
\begin{equation} \label{AlphaBetaGamma}
	\boxed{ \alpha_t = \log{p} - \log{t} , \qquad  \beta_t = p \log{t} + \log{C},  \qquad \gamma_t = p \log{t}, }
\end{equation}
where $C>0$ is a constant. This yields the $p$-Bregman Lagrangian and Hamiltonian given by
\begin{equation} 
		\boxed{ \mathcal{L}_{p}(X,V,t) = \frac{t^{\lambda^{-1} \zeta p +1}}{2p} \langle V , V\rangle  - Cpt^{(\lambda^{-1}\zeta +1)p-1} f(X), } \end{equation} 
	\begin{equation} \label{pBregmanH} \boxed{ \mathcal{H}_{p}(X,R,t)= \frac{p}{2t^{\lambda^{-1}\zeta p +1}} \llangle R , R\rrangle + Cpt^{(\lambda^{-1}\zeta +1)p-1} f(X), }
\end{equation}
and the corresponding $p$-Bregman Euler--Lagrange equations are given by
\begin{equation} \label{eq: EL Convex} \boxed {\nabla_{\dot{X}}\dot{X}  + \frac{\zeta p +\lambda}{\lambda t} \dot{X} + Cp^2t^{p-2} \text{gradf}(X) = 0. } \end{equation} 

\begin{theorem}
Suppose that $f : \mathcal{Q} \rightarrow \mathbb{R}$ is a geodesically weakly-quasi-convex function, and that Assumption~\ref{Assumption 1} is satisfied. Then, the $p$-Bregman Euler--Lagrange equation~\eqref{eq: EL Convex} has a solution, and any solution $X(t)$ converges to a minimizer $x^*$ of $f$ with rate $\boxed{f(X(t)) - f(x^*) \leq \mathcal{O}(1/t^p)}.$
\proof{ See Appendix~\ref{Section: Appendix Existence} for the existence of a solution to the $p$-Bregman Euler--Lagrange equations. The $\mathcal{O}(1/t^p)$ convergence rate follows directly from Theorem~\ref{Thm: ConvergenceRateC}. } 
\end{theorem}  
Note that this theorem reduces to Theorem 5 from \cite{alimisis2020} when $p=2$ and $C=1/4$.  \pagebreak

\begin{remark}
To construct this variational framework for accelerated optimization, we first constructed candidate $p$-equations with the desired $\mathcal{O}(1/t^p)$ convergence rates, and then designed Lagrangians whose $p$-Bregman Euler--Lagrange equations matched the candidate $p$-equations, by inspection. We then used a similar approach to extend these results to the general $\alpha,\beta,\gamma$ case presented here.
\end{remark}

\begin{remark}
	In our generalization of the Bregman Lagrangian and Hamiltonian to Riemannian manifolds, we have specialized to the case where $h(x)= \frac{1}{2} \| x \|^2$, because its Hessian $\nabla^2h(x)$ is the identity matrix, which significantly simplifies the Euler--Lagrange equations and the analysis. In addition, it avoids the complication of making intrinsic sense of terms like $X+ e^{-\alpha} V$ in the vector space Bregman Lagrangians and Hamiltonians, which requires the use of Riemannian geodesics and exponentials since $X\in \mathcal{Q}$ while $V \in T_X \mathcal{Q}$.
\end{remark}

  \hfill 

\subsection{Strongly Convex Case}

Suppose $f : \mathcal{Q} \rightarrow \mathbb{R}$ is a geodesically $\mu$-strongly-convex function, and that Assumption~\ref{Assumption 1} is satisfied. With $\zeta$ given by equation~\eqref{eq: zeta}, let
\begin{equation}
	\eta = \left(\frac{1}{\sqrt{\zeta}} + \sqrt{\zeta} \right) \sqrt{\mu}. 
\end{equation}  
We define the corresponding Lagrangian  $\mathcal{L}^{SC} :   T\mathcal{Q} \times \mathbb{R} \rightarrow  \mathbb{R}$ by
\begin{equation} \label{BregmanLGeneralSC}
	\boxed{\mathcal{L}^{SC}(X,V,t) = \frac{e^{\eta t} }{2} \langle V , V\rangle  - e^{\eta t} f(X),}
\end{equation}
\normalsize and the corresponding Hamiltonian $\mathcal{H}^{SC}:   T^*\mathcal{Q} \times \mathbb{R} \rightarrow  \mathbb{R} $ is given by
\begin{equation} \label{BregmanHGeneralSC}
	\boxed{\mathcal{H}^{SC}(X,R,t)= \frac{e^{- \eta t} }{2} \llangle R , R\rrangle + e^{\eta t} f(X).}
\end{equation}
\normalsize 

\begin{theorem} 
	The Euler--Lagrange equation corresponding to the Lagrangian $\mathcal{L}^{SC}$ is given by
	\begin{equation} \label{eq: EL SC}
		\boxed{\nabla_{\dot{X}}\dot{X}  +\eta \dot{X} + \emph{gradf}(X) = 0.} \end{equation}
	\proof{The derivation of the Euler--Lagrange equation is presented in Appendix~\ref{Section: Appendix Derivation EL SC}.}
\end{theorem}

\begin{theorem}\label{Thm: ConvergenceRateSC}
	Suppose $f : \mathcal{Q} \rightarrow \mathbb{R}$ is a geodesically $\mu$-strongly-convex function, and suppose that Assumption~\ref{Assumption 1} is satisfied. Then, the Euler--Lagrange equation~\eqref{eq: EL SC}
	has a solution, and any solution $X(t)$ converges to a minimizer $x^*$ of $f$ with rate
	\begin{equation}
		\boxed{f(X(t)) - f(x^*) \leq \frac{\mu \| \emph{Log}_{x_0}{(x^*)} \|^2 + 2\left( f(x_0) -f(x^*)\right) }{2e^{\sqrt{\frac{\mu}{\zeta}}t}}.}
	\end{equation} 
	
	\normalsize \proof{ See Appendix~\ref{Section: Appendix Existence SC} for the existence of a solution to the Euler--Lagrange equation~\eqref{eq: EL SC}, and Theorem 7 from \cite{alimisis2020} for the convergence rate. }
\end{theorem} 

\section{Numerical Experiments} \label{Section: Numerical Experiments}

The $p$-Bregman Euler--Lagrange equation~\eqref{eq: EL Convex} can be rewritten as the first-order system
\begin{align}
	\dot{X}  = V, \qquad \quad 
	\nabla_V V   = -  \frac{\zeta p + \lambda }{\lambda t} V - Cp^2t^{p-2} \text{gradf}(X),
\end{align}
for the geodesically $\lambda$-weakly-quasi-convex case, and the Euler--Lagrange equation~\eqref{eq: EL SC} corresponding to the Lagrangian $\mathcal{L}^{SC}$ can be rewritten as the first-order system
\begin{align}
	\dot{X}  = V, \qquad \quad 
	\nabla_V V   = -  \left(\frac{1}{\sqrt{\zeta}} + \sqrt{\zeta} \right) \sqrt{\mu}  V - \text{gradf}(X),
\end{align}
for the $\mu$-strongly convex case. As in \cite{alimisis2020}, we can adapt a semi-implicit Euler scheme (explicit Euler update for the velocity $V$ followed by an update for position $X$ based on the updated value of $V$) to the Riemannian setting to obtain the following algorithm: 

\begin{algorithm}[H]
	\DontPrintSemicolon
	
	\KwInput{A function $f : \mathcal{Q} \rightarrow \mathbb{R}$. Constants $C,h,p>0$. $X_0 \in \mathcal{Q}$. $V_0 \in T_{X_0} \mathcal{Q}$. }
	
		\While{convergence criterion is not met}
	{
		\uIf{$f$ is $\mu$-geodesically strongly convex}
		{
			$b_k \leftarrow  1- h  \left(\frac{1}{\sqrt{\zeta}} + \sqrt{\zeta} \right) \sqrt{\mu}, \quad  c_k \leftarrow 1$\;
}	

\ElseIf{$f$ is $\lambda$-weakly-quasi-convex}{
		$b_k  \leftarrow 1 -  \frac{\zeta p+\lambda}{ \lambda k}, \quad c_k \leftarrow C p^2 (kh)^{p-2} $\;
	}
		\textbf{Version I}:	$a_k \leftarrow b_k V_k - hc_k \text{gradf}(X_k)$ \;
		\textbf{Version II}:	$a_k \leftarrow b_k V_k - h c_k \text{gradf}\left(\text{Exp}_{X_k}(hb_k V_k)\right)$ \;
		
		$X_{k+1} \leftarrow \text{Exp}_{X_k}(ha_k), \quad V_{k+1} \leftarrow \Gamma_{X_k}^{X_{k+1}} a_k$ \;
	} 
	\caption{Semi-Implicit Euler Integration of the $p$-Bregman Euler--Lagrange Equations } \label{Alg:Semi-Implicit_Euler}
\end{algorithm} 

Version I of Algorithm~\ref{Alg:Semi-Implicit_Euler} corresponds to the usual update for the Semi-Implicit Euler scheme, while Version II is inspired by the reformulation of Nesterov's method from \cite{Sutskever2013} that uses a corrected gradient $\nabla f(X_k +h b_kV_k)$ instead of the traditional gradient $\nabla f(X_k)$. Note that the SIRNAG algorithm presented in \cite{alimisis2020} corresponds to the special case where $p=2$ and $C=1/4$. 

The first problem we have investigated is the problem presented in \cite{alimisis2020} of minimizing the (strongly convex) distance function $f(x) =\frac{1}{2} d(x,q)^2$ for a given point $q$, on a subset of chosen finite diameter of the hyperbolic plane $\mathbb{H}^2$, which is a manifold with constant negative curvature $K=-1$. 

The second problem we have investigated is Rayleigh quotient optimization. Eigenvectors corresponding to the largest eigenvalue of a symmetric $n \times n$ matrix $A$ maximize the Rayleigh quotient $\frac{v^\top Av}{ v^\top v}$ over $\mathbb{R}^n$. Thus, a unit eigenvector $v^*$ corresponding to the largest eigenvalue of the matrix $A$ is a minimizer of the function $f(v) = -  v^\top Av,$ over the unit sphere $\mathcal{Q} = \mathbb{S}^{n-1}$, which can be thought of as a Riemannian submanifold with constant positive curvature $K=1$ of $\mathbb{R}^{n}$ endowed with the Riemannian metric inherited from the Euclidean inner product $g_v(u,w) = u^\top w$. More information concerning the geometry of $\mathbb{S}^{n-1}$, such as its tangent bundle, its orthogonal projection and exponential map can be found in  \cite{Absil2008}. Solving the Rayleigh quotient optimization problem efficiently is challenging when the given symmetric matrix $A$ is ill-conditioned and high-dimensional. Note that an efficient algorithm that solves the above minimization problem can also be used to find eigenvectors corresponding to the smallest eigenvalue of $A$ by using the fact that the eigenvalues of $A$ are the negative of the eigenvalues of $-A$. 

\begin{figure}[!ht] 
	\centering
	\begin{minipage}[b]{0.49\textwidth}
		\includegraphics[width=\textwidth]{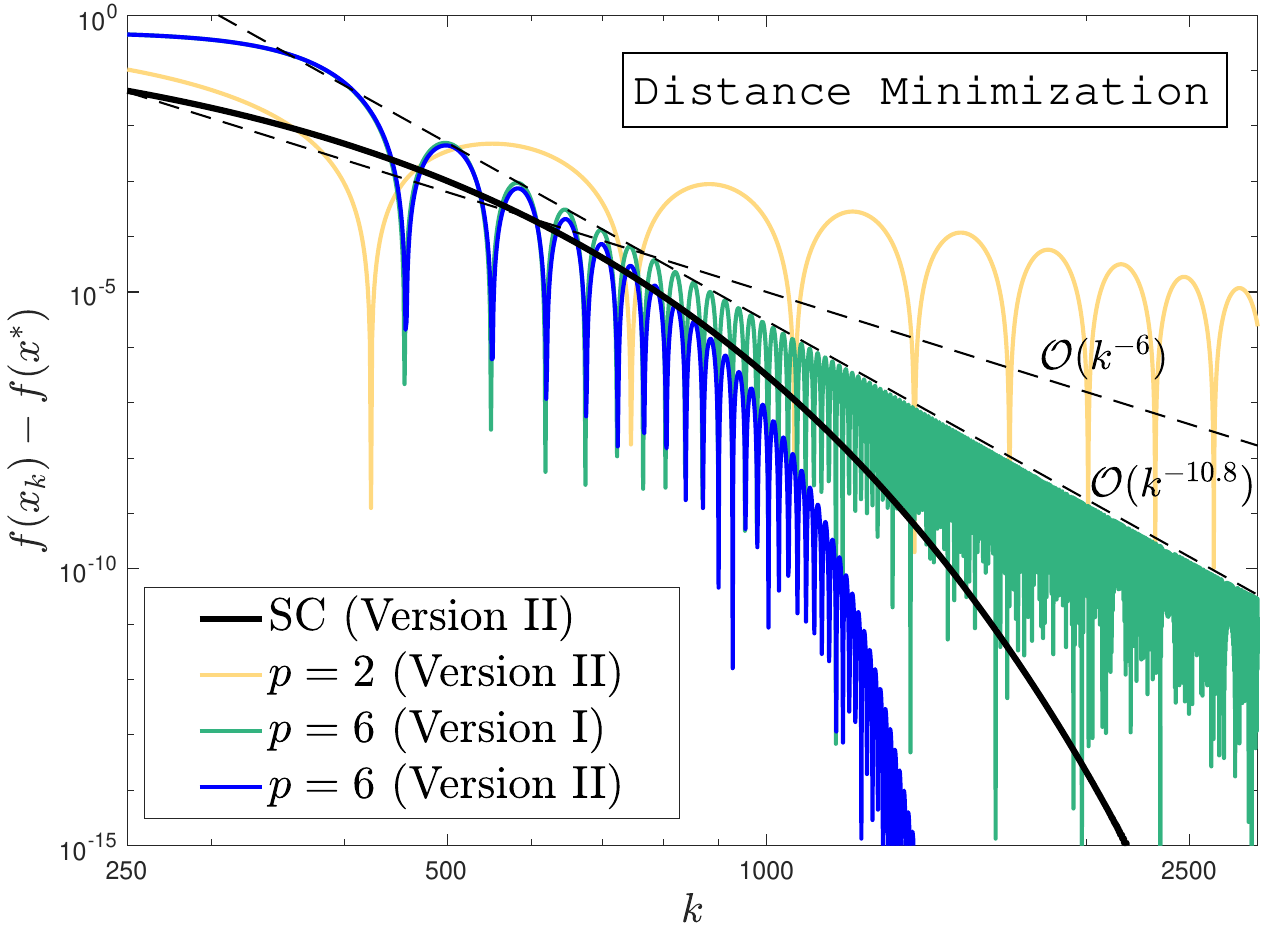}
	\end{minipage}
	\begin{minipage}[b]{0.49\textwidth}
		\includegraphics[width=\textwidth]{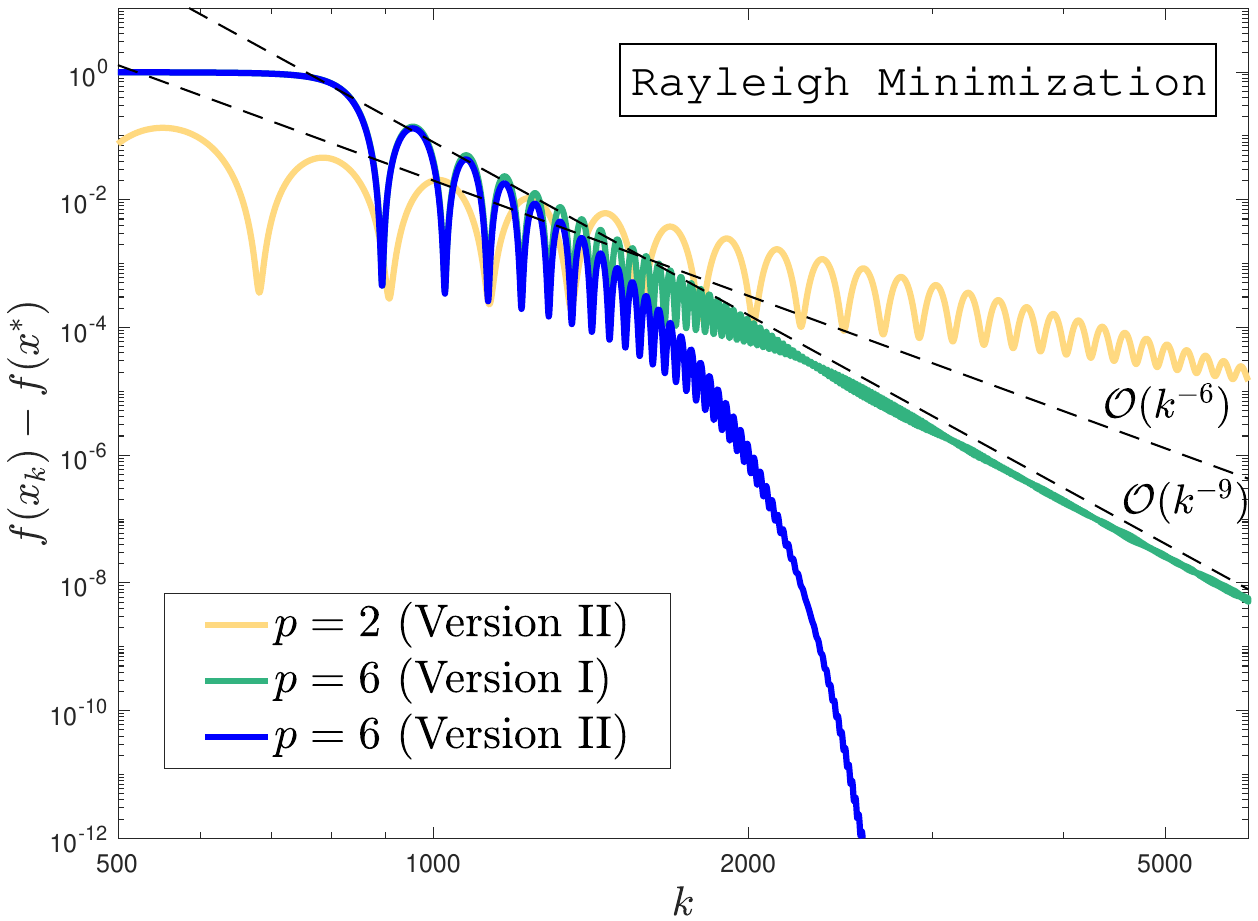}
	\end{minipage}
	\caption{Comparison of the rates of convergence of the $\mu$-strongly convex (SC) Algorithm~\ref{Alg:Semi-Implicit_Euler} and convex Algorithms~\ref{Alg:Semi-Implicit_Euler} with different values of $p$ and with the two versions of the update corresponding to the traditional and corrected gradients. Note that all the algorithms were implemented with the same time-step $h$.
	}\label{fig: ConvergenceRates}
\end{figure}

Experiments carried out in \cite{alimisis2020} showed that SIRNAG (the convex $p=2$ Algorithm~\ref{Alg:Semi-Implicit_Euler}) and the strongly convex Algorithm~\ref{Alg:Semi-Implicit_Euler} were of comparable efficiency or more efficient than the standard Riemannian Gradient Descent (RGD) method, depending on the properties of the objective function and on the geometry of the Riemannian manifold. We have conducted further numerical experiments to investigate how the simple discretization of higher-order $p=6$ Bregman dynamics compared to its $p=2$ counterpart, and to see whether it matches the $\mathcal{O}(k^{-p})$ convergence rate. The numerical results obtained for the distance minimization and Rayleigh minimization problems are illustrated in Figure~\ref{fig: ConvergenceRates}, where all the algorithms were implemented with the same fixed time-step. We can see that the $p=6$ algorithms outperform their $p=2$ counterparts, and that the efficiency improvement is very important. Furthermore, both versions of the $p=6$ Algorithm~\ref{Alg:Semi-Implicit_Euler} exhibit a faster convergence rate than $\mathcal{O}(k^{-6})$. While Version I of Algorithm~\ref{Alg:Semi-Implicit_Euler} exhibits polynomial rates of $\mathcal{O}(k^{-10.8}) $ and $\mathcal{O}(k^{-9}) $ on the objective functions considered, Version II of Algorithm~\ref{Alg:Semi-Implicit_Euler} exhibits a much faster exponential rate of convergence on both examples. 

\begin{figure}[!ht] 
	\centering
	\begin{minipage}[b]{0.94\textwidth}
		\includegraphics[width=\textwidth]{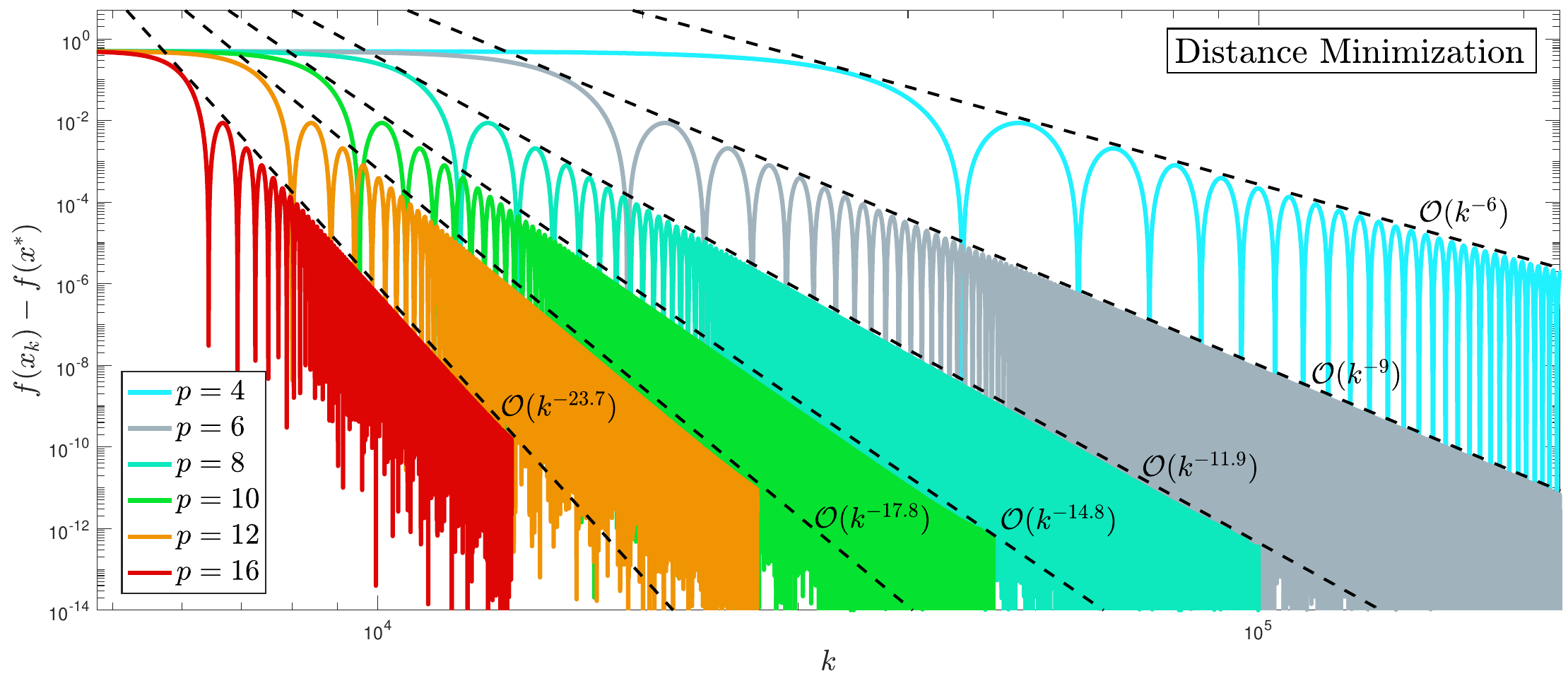}
	\end{minipage}
	\begin{minipage}[b]{0.94\textwidth}
		\includegraphics[width=\textwidth]{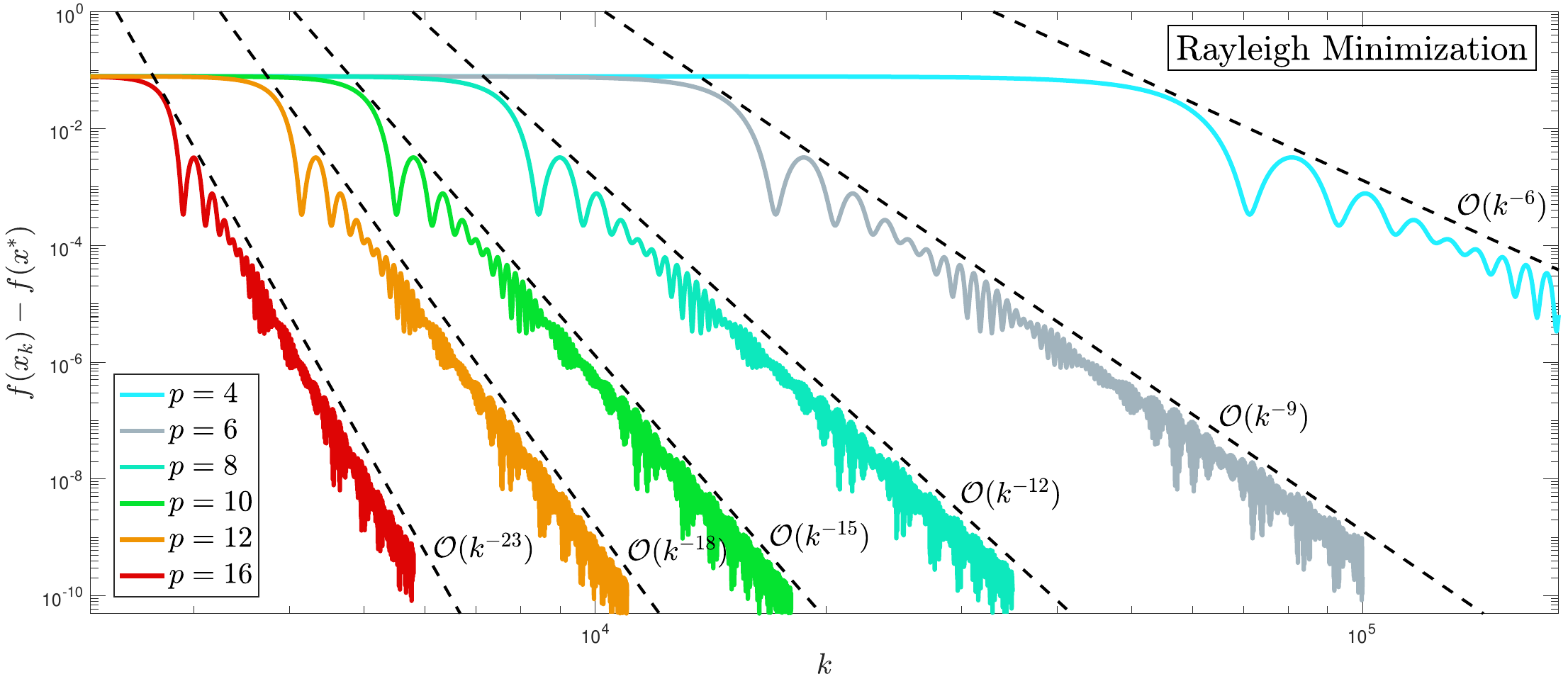}
	\end{minipage}
	\caption{Evolution of the rates of convergence of Version 1 of the convex Algorithm~\ref{Alg:Semi-Implicit_Euler} with different values of $p$. Note that all the algorithms were implemented with the same time-step $h$.
	}\label{fig: Evolution}
\end{figure}

Figure~\ref{fig: Evolution} displays the evolution of the rates of convergence of Version 1 of the convex Algorithm~\ref{Alg:Semi-Implicit_Euler} as the value of the parameter $p$ is increased from $p=4$ to $p=16$ for the distance minimization and Rayleigh minimization problems. We can clearly see an improvement in the convergence rates as the value of $p$ increases, and for each value of $p$ the algorithm achieves a faster rate of convergence than $\mathcal{O}(k^{-p})$.

Note however that an increase in the value of $p$ in Algorithm~\ref{Alg:Semi-Implicit_Euler}, which corresponds to an increase in the order of the Bregman dynamics integrated, requires a decrease in the time-step, in agreement with intuitive expectations. This time-step decrease requirement is especially important due to the polynomially growing $h(kh)^{p-2}$ coefficient multiplying the gradient of $f$ in the updates of the algorithm. Such a decrease in the time-step does not really affect the convergence rate, but the transition between the initialization and convergence phases takes longer. As a consequence, by using larger time-steps, the algorithm corresponding to a smaller value of $p$ might achieve a desired convergence criterion with fewer iterations than the algorithm corresponding to a larger value of $p$, despite having a slower convergence rate. Similar issues arise when discretizing the continuous Euler--Lagrange flow associated with accelerated optimization on vector spaces, and in that situation, it was observed that time-adaptive symplectic integrators based on Hamiltonian variational integrators resulted in dramatically improved robustness and stability. As such, it will be natural to explore generalizations of time-adaptive symplectic integrators based on Hamiltonian variational integrators applied to Poincar\'e transformed Hamiltonians, that respect the Riemannian manifold structure in order to yield more robust and stable numerical discretizations of the flows we have studied in this paper in order to construct accelerated optimization algorithms on Riemannian manifolds. We will lay the foundation for such time-adaptive symplectic integrators in Section~\ref{sec: TimeInvariance}.

\begin{figure}[!ht] 
	\centering
	\begin{minipage}[b]{0.49\textwidth}
		\includegraphics[width=\textwidth]{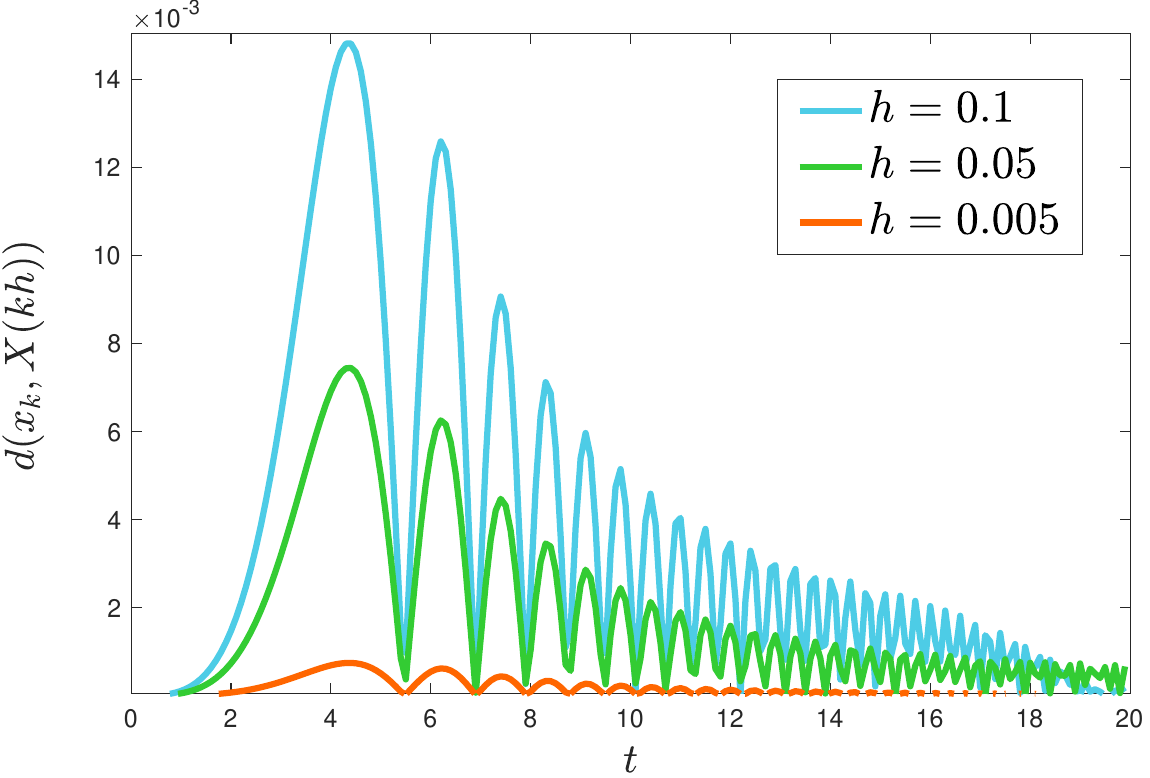}
	\end{minipage}
	\begin{minipage}[b]{0.49\textwidth}
		\includegraphics[width=\textwidth]{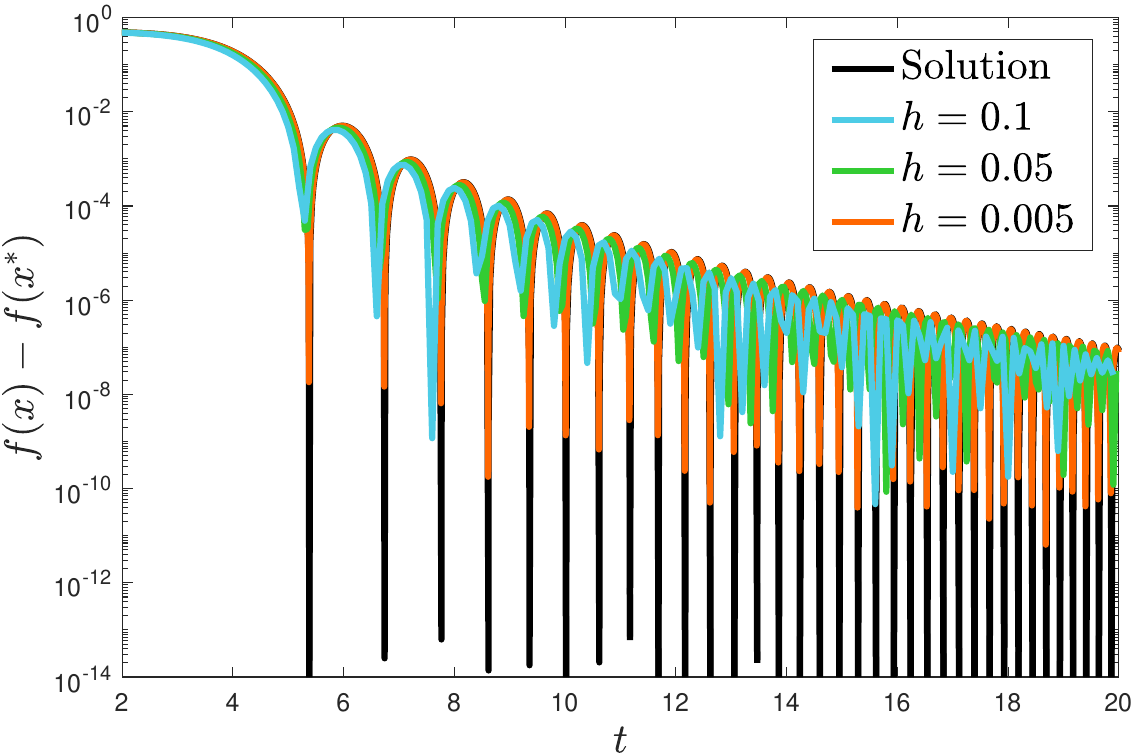}
	\end{minipage}
	\caption{Discretization errors (top graph) and convergence rates (bottom graphs) of Version I of the $p=5$ convex Algorithm~\ref{Alg:Semi-Implicit_Euler} with different values of $h$ for the distance minimization problem. The true solution of the differential equation was approximated by the same algorithm with a very small time-step $h=10^{-5}$.	}\label{fig: Discretizations}
\end{figure}

Finally, Figure~\ref{fig: Discretizations} shows that the discretization empirically converges to the solution of the ODE as the time-step $h$ goes to 0. Note that although all the discretizations follow the ODE trajectory closely, smaller time-steps result in a larger number of iterations, especially to transition from the initialization plateau to the convergence phase (around time $t=4$ in the example presented in Figure~\ref{fig: Discretizations}). A theoretical shadowing result bounding the error between the discrete-time RGD and its continuous-time limiting ODE was obtained in \cite{alimisis2020} thanks to the uniform contraction property of the dynamical system associated with Riemannian Gradient Descent. It would be desirable to obtain similar shadowing results in the future for discretizations of the class of ODEs considered in this paper, perhaps drawing inspiration from \cite{Jadbabaie2018}. However, such a result might be very difficult to obtain because momentum methods lack contraction, are nondescending, and are highly oscillatory~\cite{alimisis2020,Orvieto2019}. While it is hoped that the continuous analysis in this paper will eventually guide the convergence analysis of discrete-time algorithms, this does not appear to be a straightforward exercise, as one would first need to reconcile the arbitrarily fast $\mathcal{O}(1/t^p)$ rate of convergence of the continuous-time trajectories with Nesterov's barrier theorem of $\mathcal{O}(1/k^2)$ for discrete-time algorithms. Even on normed vector spaces, obtaining theoretical guarantees was a challenging task, achieved in \cite{Jadbabaie2018} in the special case where $p>2$ under additional assumptions on the objective function and on its derivatives. Generalizing these results to the general family of $\alpha,\beta,\gamma$ Bregman Lagrangians on Riemannian manifolds would be much more challenging since the notions of derivatives become more complicated, and since all the usual vector space operations and objects have to be replaced by their Riemannian generalization which involve geodesics, parallel transport, Riemannian exponentials and Riemannian logarithms. 

\section{Time Invariance and Poincar\'e Transformation} \label{sec: TimeInvariance}

Let $f : \mathcal{Q} \rightarrow \mathbb{R}$ be a given $\lambda$-weakly-quasi-convex function, and suppose Assumption~\ref{Assumption 1} is satisfied. In Section~\ref{Sec: VariationalFormulation}, we formulated a variational framework for the minimization of $f$, via Bregman Lagrangians and Hamiltonians. We now extend Theorem~\ref{ThmTimeDilation} to Riemannian manifolds. 
\begin{theorem} 
	Suppose that Assumption~\ref{Assumption 1} is satisfied and that the curve $X(t)$ satisfies the Riemannian Bregman Euler--Lagrange equation~\eqref{eq: EL} corresponding to $\mathcal{L}_{\alpha,\beta,\gamma}$. Then the reparametrized curve $ X(\tau(t))$ satisfies the Bregman Euler--Lagrange equation~\eqref{eq: EL} corresponding to the modified Riemannian Bregman Lagrangian $\mathcal{L}_{\tilde{\alpha},\tilde{\beta},\tilde{\gamma}}$ where $\tilde{\alpha}_t = \alpha_{\tau(t)} + \log{\dot{\tau}(t) }$,  $\tilde{\beta}_t = \beta_{\tau(t)}$, and $\tilde{\gamma}_t = \gamma_{\tau(t)}$. Furthermore $\alpha,\beta,\gamma$ satisfy the ideal scaling conditions (\ref{IdealScaling}) if and only if $\tilde{\alpha},\tilde{\beta},\tilde{\gamma}$ do.
	\proof{See Appendix~\ref{Section: Appendix Proof Invariance}.}
\end{theorem} 

As a special case, we have the following theorem:

\begin{theorem} 
	Suppose that $f : \mathcal{Q} \rightarrow \mathbb{R}$ is a geodesically $\lambda$-weakly-quasi-convex function, and that Assumption~\ref{Assumption 1} is satisfied.  Suppose $X(t)$ satisfies the $p$-Bregman Euler--Lagrange equation~\eqref{eq: EL Convex}.
	Then, the reparametrized curve $X(t^{\mathring{p}/p})$ satisfies the $\mathring{p}$-Bregman Euler--Lagrange equation~\eqref{eq: EL Convex}.
\end{theorem}

Thus, the entire subfamily of Bregman trajectories indexed by the parameter $p$ can be obtained by speeding up or slowing down along the Bregman curve in spacetime corresponding to any specific value of $p$. Inspired by the computational efficiency of the approach introduced in \cite{duruisseaux2020adaptive}, it is natural to attempt to exploit the time-rescaling property of the Bregman dynamics together with a carefully chosen Poincar\'e transformation to transform the $p$-Bregman Hamiltonian into an autonomous version of the $\mathring{p}$-Bregman Hamiltonian in extended phase-space, where $\mathring{p} < p$. This would allow us to integrate the higher-order $p$-Bregman dynamics while benefiting from the computational efficiency of integrating the lower-order $\mathring{p}$-Bregman dynamics. Explicitly, the time rescaling $\tau(t) = t^{\mathring{p}/p}$ is associated to the monitor function 
\begin{equation} \label{eq: TimeTransformation}
	\frac{dt}{d\tau} = g_{p \rightarrow \mathring{p}}(t) = \frac{p}{\mathring{p}} t^{1-\mathring{p}/p},
\end{equation}
and generates a Poincar\'e transformed Hamiltonian
\begin{equation} \label{eq: NewHamiltonian}
	\bar{\mathcal{H}}_{p \rightarrow \mathring{p}}(\bar{X},\bar{R}) = g_{p \rightarrow \mathring{p}}(X^t) \left(\mathcal{H}_p\left(\bar{X},R\right) + R^t \right),
\end{equation}
in the extended space $\bar{\mathcal{Q}} = \mathcal{Q}\times \mathbb{R}$ where $\bar{X} = \begin{bmatrix} X \\ X^t \end{bmatrix} $ and $\bar{R} = \begin{bmatrix} R \\ R^t \end{bmatrix} $. We will make the conventional choice $X^t=t$, with conjugate momentum $R^t$, and $R^t(0)=-\mathcal{H}_p(X(0),R(0),0) = -H_0$, which is chosen so that $	\bar{\mathcal{H}}_{p \rightarrow \mathring{p}}(\bar{X},\bar{R}) =0$ along all integral curves through $(\bar{X}(0),\bar{R}(0))$. The time $t$ shall be referred to as the physical time, while $\tau$ will be referred to as the fictive time. The corresponding Hamiltonian equations of motion in the extended phase space are then given by
\begin{align}
	\dot{\bar{X}}  = \frac{\partial \bar{\mathcal{H}}_{p \rightarrow \mathring{p}}}{\partial \bar{R}},  \qquad\quad  \dot{\bar{R}} =  -\frac{\partial \bar{\mathcal{H}}_{p \rightarrow \mathring{p}}}{\partial \bar{X}}.    
\end{align}
Now, suppose $(\bar{X}(\tau) , \bar{R}(\tau))$ are solutions to these extended equations of motion, and let $(x(t),r(t))$ solve Hamilton's equations for the original Hamiltonian $\mathcal{H}_p$. Then
\begin{equation*}
\bar{\mathcal{H}}_{p \rightarrow \mathring{p}}(\bar{X}(\tau) , \bar{R}(\tau))  = \bar{\mathcal{H}}_{p \rightarrow \mathring{p}}(\bar{X}(0) , \bar{R}(0)) = 0.\end{equation*}
Thus, the components $(X(\tau) , R(\tau))$ in the original phase space of $(\bar{X}(\tau) , \bar{R}(\tau))$ satisfy
\[\mathcal{H}_p(X(\tau) , R(\tau) ,\tau ) = -R^t(\tau), \qquad  \mathcal{H}_p(X(0) , R(0) ,0 ) = -R^t(0) = \mathcal{H}_p(x(0),r(0),0). \]
Therefore,  $(X(\tau) , R(\tau))$ and $(x(t),r(t))$ both satisfy Hamilton's equations for the original Hamiltonian $\mathcal{H}_p$ with the same initial values, so they must be the same. 

As a consequence, instead of integrating the $p$-Bregman Hamiltonian system~\eqref{pBregmanH}, we can focus on the Poincar\'e transformed Hamiltonian $\bar{\mathcal{H}}_{p \rightarrow \mathring{p}}$ in extended phase-space given by equation~\eqref{eq: NewHamiltonian}, with $\mathcal{H}_p$ and $g_{p \rightarrow \mathring{p}}$ given by equations~\eqref{pBregmanH} and \eqref{eq: TimeTransformation}, that is
\begin{equation}  \label{eq: Transformed H}
	\boxed{\bar{\mathcal{H}}_{p \rightarrow \mathring{p}}(\bar{X},\bar{R})  = 	 \frac{p^2}{2\mathring{p} (X^t)^{\lambda^{-1} \zeta p +\mathring{p}/p}}  \llangle R , R\rrangle  + \frac{Cp^2}{\mathring{p}}(X^t)^{(\lambda^{-1} \zeta  +1)p-\mathring{p}/p} f(X)+ \frac{p}{\mathring{p}} (X^t)^{1-\mathring{p}/p}   R^t, }
\end{equation}
The resulting
integrator has constant time-step in fictive time $\tau$ but variable time-step in physical time $t$. In our prior work on discretizations of variational formulations of accelerated optimization on normed spaces~\cite{duruisseaux2020adaptive}, we performed a very careful computational study of how time-adaptivity and symplecticity of the numerical scheme improve the performance of the resulting numerical optimization algorithm. In particular, we observed that time-adaptive Hamiltonian variational discretizations, which are automatically symplectic, with adaptive time-steps informed by the time invariance of the family of $p$-Bregman Lagrangians and Hamiltonians yielded the most robust and computationally efficient numerical optimization algorithms, outperforming fixed-timestep symplectic discretizations, adaptive-timestep non-symplectic discretizations, and Nesterov's accelerated gradient algorithm which is neither time-adaptive nor symplectic. As such, it would be desirable to generalize the time-adaptive Hamiltonian variational integrator framework to Riemannian manifolds, and apply it to the variational formulation of accelerated optimization on Riemannian manifolds.

Note that the variational framework for accelerated optimization presented in Section~\ref{Sec: VariationalFormulation} has also been exploited successfully in the special case of Lie groups in subsequent papers~\cite{Lee2021,Duruisseaux2022Lagrangian}, using two different formulations of time-adaptive symplectic Lagrangian integration, with very promising numerical results. Another important case involves Riemannian submanifolds that are embedded in a Riemannian linear manifold and are realized as the level set of a submersion. The characterization of the submanifold as the level set of a submersion, together with the linear space structure of the embedding space, and the variational characterization of the dynamics naturally lends itself to the use of the Lagrange multiplier theorem, which allows one to use Hamiltonian variational integrators defined on the embedding space by including a Lagrange multiplier term involving the submersion in the Lagrangian or Hamiltonian \cite{Duruisseaux2021Constrained}. This is analogous to the derivation of the SHAKE and RATTLE methods as variational integrators for constrained systems (see, for example, Section~3.5 of~\cite{MaWe2001}). Another practical method can be obtained by projecting the updates of Hamiltonian variational integrators defined on the embedding space onto the constraint manifold~\cite{Duruisseaux2021Projection}. The numerical results in these subsequent papers \cite{Duruisseaux2021Constrained,Duruisseaux2021Projection} suggest that the time-adaptive Hamiltonian approach can be very competitive when numerically solving optimization problems on Riemannian manifolds.

\section{Conclusion}

We have shown that on Riemannian manifolds, the convergence rate in continuous time of a geodesically convex or weakly-quasi-convex function $f(x(t))$ to its optimal value can be accelerated to an arbitrary convergence rate, which extended the results of \cite{WiWiJo16} from normed vector spaces to Riemannian manifolds. This rate of convergence is achieved along solutions of the Euler--Lagrange and Hamilton's equations corresponding to a family of time-dependent Bregman Lagrangian and Hamiltonian systems on Riemannian manifolds. As was demonstrated in the normed vector space setting, such families of Bregman Lagrangians and Hamiltonians can be used to construct practical, robust, and computationally efficient numerical optimization algorithms that outperform Nesterov's accelerated gradient method by considering geometric structure-preserving discretizations of the continuous-time flows.

Numerical experiments implementing a simple discretization of the $p$-Bregman Euler--Lagrange equations applied to a distance minimization and Rayleigh minimization problems confirmed that the higher-order algorithms outperform significantly their lower-order counterparts and the corresponding $\mathcal{O}(1/k^{p})$ convergence rates. Numerical results also showed that using a corrected gradient in the update instead of the traditional gradient, as was done in \cite{Sutskever2013}, improved the theoretically predicted polynomial convergence rate to an exponential rate of convergence in practice. While higher values of $p$ result in faster rates of convergence, they usually require smaller time-steps and also appear to be more prone to stability issues under numerical discretization, which can cause the numerical optimization algorithm to diverge, but we anticipate that symplectic discretizations will address these stability issues. 

Finally, in analogy to what was done in \cite{WiWiJo16} for normed vector spaces, we proved that the family of time-dependent Bregman Lagrangian and Hamiltonians on Riemannian manifolds is closed under time rescaling. Inspired by the computational efficiency of the approach introduced in \cite{duruisseaux2020adaptive}, we can then exploit this invariance property via a carefully chosen Poincar\'e transformation that will allow us to integrate higher-order $p$-Bregman dynamics while benefiting from the computational efficiency of integrating a lower-order $\mathring{p}$-Bregman Hamiltonian system. 

It was observed in our prior computational experiments in the normed vector space case~\cite{duruisseaux2020adaptive} that geometric discretizations which respect the time-rescaling invariance and symplecticity of the Bregman Lagrangian and Hamiltonian flows were substantially less prone to stability issues, and were therefore more robust, reliable, and computationally efficient. As such, it is natural to develop time-adaptive Hamiltonian variational integrators for the Bregman Hamiltonian introduced in this paper describing accelerated optimization on Riemannian manifolds.

Developing an intrinsic extension of Hamiltonian variational integrators to manifolds will require some additional work, since the current approach involves Type II/Type III generating functions $H_d^+(q_k, p_{k+1})$, $H_d^-(p_k, q_{k+1})$, which depend on the position at one boundary point, and the momentum at the other boundary point. However, this does not make intrinsic sense on a manifold, since one needs the base point in order to specify the corresponding cotangent space, and one should ideally consider a Hamiltonian variational integrator construction based on discrete Dirac mechanics~\cite{LeOh2008}, which would yield a generating function $E_d^+(q_k, q_{k+1}, p_{k+1})$, $E_d^-(q_k, p_k, q_{k+1})$, that depends on the position at both boundary points and the momentum at one of the boundary points. This approach can be viewed as a discretization of the generalized energy $E(q,v,p)=\langle p,v\rangle - L(q,v)$, in contrast to the Hamiltonian $H(q,p)=\ext_{v}\langle p,v\rangle - L(q,v)=\left.\langle p,v\rangle - L(q,v)\right|_{p=\frac{\partial L}{\partial v}}$.

However, a more practical method relies on the fact that we have a Riemannian manifold, which is endowed with a Riemannian exponential and Riemannian logarithm that can be used to construct an extension of Hamiltonian variational integrators using geodesic normal coordinates. For many important matrix manifolds, one can replace the Riemannian exponential in the geodesic normal coordinates by a retraction~\cite{Absil2008}, which is often constructed using matrix factorizations.

We anticipate that applying an appropriate generalization of Hamiltonian variational integrators to the Bregman Hamiltonians introduced in this paper will yield a novel class of robust and efficient accelerated optimization algorithms on Riemannian manifolds. The variational framework for accelerated optimization presented in Section~\ref{Sec: VariationalFormulation} has also been exploited successfully in the special case of Lie groups in subsequent papers \cite{Lee2021,Duruisseaux2022Lagrangian}, using two different formulations of time-adaptive symplectic Lagrangian integration, with very promising numerical results which illustrate that our framework can be very competitive for optimization problems of interest on Lie groups and more generally on Riemannian manifolds. As mentioned at the end of Section~\ref{sec: TimeInvariance}, another important case involves Riemannian submanifolds that are embedded in a Riemannian linear manifold and are realized as the level set of a submersion. In~\cite{Duruisseaux2021Constrained}, we studied how holonomic constraints can be incorporated into variational integrators to constrain the updates of the numerical optimization algorithm to the Riemannian manifold of interest, and in \cite{Duruisseaux2021Projection}, the manifold constraints were enforced via projections. The numerical results in these two subsequent papers suggest that the time-adaptive Hamiltonian approach introduced in this paper can be the basis for competitive numerical optimization algorithms on Riemannian manifolds.

It would be desirable in future work to analyze the resulting discrete-time algorithms and rigorously establish their rates of convergence. Although theoretical shadowing results have already been derived for certain discrete optimization algorithms on Riemannian manifolds, such a result might be very difficult to obtain for the momentum-based algorithms presented in this paper because momentum methods lack contraction, are nondescending and highly oscillatory~\cite{alimisis2020,Orvieto2019}. It might also be possible to generalize the theoretical guarantees obtained laboriously on normed vector spaces in~\cite{Jadbabaie2018}, but this would be an even more challenging task since the usual vector space operations and objects have to be replaced by their more convoluted Riemannian generalizations. In addition, we would like to better understand how to reconcile the arbitrarily high rate of convergence one expects from the continuous-time analysis, with Nesterov's barrier theorem on the rate of convergence of discrete-time algorithms.

\section*{Acknowledgments} 

The authors would like to thank the referees for their careful review of this paper and their helpful suggestions. \\

\appendix

\section{Derivation of the Euler--Lagrange Equations}

\subsection{Convex and Weakly-Quasi-Convex Cases}\label{Section: Appendix Derivation EL}

\begin{theorem}
	The Euler--Lagrange equation corresponding to the Lagrangian
	\[
	\mathcal{L}_{\alpha,\beta,\gamma}(X,V,t) = \frac{1}{2}e^{ \lambda^{-1}\zeta \gamma_t - \alpha_t} \langle V , V\rangle  - e^{\alpha_t + \beta_t +  \lambda^{-1}\zeta \gamma_t} f(X),
	\]
	is given by
	\[
	\nabla_{\dot{X}}\dot{X}  +\left( \lambda^{-1} \zeta e^{\alpha_t} -  \dot{\alpha}_t \right) \dot{X} + e^{2\alpha_t+\beta_t }  \emph{gradf}(X) = 0,
	\] 
	\proof{ Consider a path on the manifold $\mathcal{Q}$ described in coordinates by
		\[ \left( x(t),\dot{x}(t) \right) = \left(  q^1(t), \ldots , q^n(t) , v^1(t), \ldots , v^n(t)   \right).\]
		Then, with $\langle \cdot , \cdot \rangle = \sum_{i,j=1}^{n}{g_{ij} dx^i dx^j},$ the Bregman Lagrangian $\mathcal{L}_{\alpha,\beta,\gamma}$ can be written as
		\[ \mathcal{L}_{\alpha,\beta,\gamma}\left( x(t),\dot{x}(t) ,t \right)  = \frac{1}{2} e^{\lambda^{-1} \zeta \gamma_t - \alpha_t} \sum_{i,j=1}^{n}{g_{ij}(x(t)) v^i(t) v^j(t) }  - e^{\alpha_t +\beta_t +\lambda^{-1} \zeta \gamma_t}f(x(t)). \]
		For $k=1,\ldots n$,
		\begin{align*} \frac{d}{dt} \left( \frac{\partial \mathcal{L}_{\alpha,\beta,\gamma} }{\partial v^k} \left( x(t),\dot{x}(t) ,t \right) \right) & =   e^{\lambda^{-1} \zeta \gamma_t - \alpha_t} \sum_{i=1}^{n}{ g_{ik}(x(t)) \frac{dv^i}{dt}(t)}   +  e^{\lambda^{-1} \zeta \gamma_t - \alpha_t} \sum_{i,j=1}^{n}{ \frac{\partial g_{kj}}{\partial q^i}(x(t)) v^i(t) v^j(t)  } \\ & \qquad \qquad \quad\qquad +  (\lambda^{-1} \zeta \dot{\gamma}_t - \dot{\alpha}_t)e^{\lambda^{-1} \zeta \gamma_t - \alpha_t} \sum_{i=1}^{n}{g_{ik}(x(t)) v^i(t)}, 
		\end{align*}
		\[  \frac{\partial \mathcal{L}_{\alpha,\beta,\gamma} }{\partial q^k}  \left( x(t),\dot{x}(t) ,t \right)  = \frac{1}{2}e^{\lambda^{-1} \zeta\gamma_t - \alpha_t} \sum_{i,j=1}^{n}{ \frac{\partial g_{ij}}{\partial q^k}(x(t)) v^i(t)v^j(t)} - e^{\alpha_t +\beta_t +\lambda^{-1} \zeta \gamma_t} \frac{\partial f}{\partial q^k} (x(t)). \]
		Multiplying both terms by $e^{\alpha_t -\lambda^{-1} \zeta \gamma_t}$, the Euler--Lagrange equations~\eqref{eq: EL Basic} for the Bregman Lagrangian $\mathcal{L}_{\alpha,\beta,\gamma}$ are given, for $k=1, \ldots , n$, by
		\begin{align*}  0 = \sum_{i=1}^{n}{ g_{ik}(x(t)) \frac{dv^i}{dt}(t)}  & +   \sum_{i,j=1}^{n}{ \frac{\partial g_{kj}}{\partial q^i}(x(t)) v^i(t) v^j(t)  } +   (\lambda^{-1} \zeta \dot{\gamma}_t - \dot{\alpha}_t) \sum_{i=1}^{n}{g_{ik}(x(t)) v^i(t)} \\ &  \qquad\quad -  \frac{1}{2} \sum_{i,j=1}^{n}{ \frac{\partial g_{ij}}{\partial q^k}(x(t)) v^i(t)v^j(t)} + e^{2\alpha_t + \beta_t} \frac{\partial f}{\partial q^k} (x(t)) .
		\end{align*}
		Rearranging terms, and multiplying by the matrix $(g^{ij})$ which is the inverse of $(g_{ij})$, we get, for $k=1,\ldots n$,  the equation
		\[ \left( \frac{dv^k}{dt}(t)   +   \sum_{i,j=1}^{n}{ \Gamma^{k}_{ij}(x(t)) v^i(t) v^j(t)  }  \right)  + \left(\lambda^{-1} \zeta \dot{\gamma}_t - \dot{\alpha}_t\right) v^k(t) + e^{2\alpha_t + \beta_t}  \left( \text{gradf}(x(t)) \right)^k  = 0,\]
		where $\Gamma^{k}_{ij}$ are the Christoffel symbols given by
		$ \Gamma^{k}_{ij} = \frac{1}{2} \sum_{l=1}^{n}{ g^{kl} \left[ \frac{\partial g_{jl}}{\partial x^i} + \frac{\partial g_{li}}{\partial x^j} - \frac{\partial g_{ij}}{\partial x^l}\right]},$ which gives the desired Euler--Lagrange equation once we use the ideal scaling equation $\dot{\gamma}_t = e^{\alpha_t}$. \qed }
\end{theorem}

\subsection{Strongly Convex Case}\label{Section: Appendix Derivation EL SC}

\begin{theorem} 
	The Euler--Lagrange equation corresponding to the Lagrangian $\mathcal{L}^{SC}$ is given by
	\begin{equation*}
		\nabla_{\dot{X}}\dot{X}  +\eta \dot{X} + \emph{gradf}(X) = 0.
	\end{equation*}	 
	\proof{Consider a path on the manifold $\mathcal{Q}$ described in coordinates by
	\[ \left( x(t),\dot{x}(t) \right) = \left(  q^1(t), \ldots , q^n(t) , v^1(t), \ldots , v^n(t)   \right).\]
	Then, with $\langle \cdot , \cdot \rangle = \sum_{i,j=1}^{n}{g_{ij} dx^i dx^j},$ the Lagrangian $\mathcal{L}^{SC}$ can be written as
	\[ \mathcal{L}^{SC}\left( x(t),\dot{x}(t) ,t \right)  = \frac{e^{\eta t}}{2} \sum_{i,j=1}^{n}{g_{ij}(x(t)) v^i(t) v^j(t) }  -  e^{\eta t} f(x(t)). \]
	For $k=1,\ldots n$, 
	\begin{align*} \frac{d}{dt} \left( \frac{\partial \mathcal{L}^{SC} }{\partial v^k} \left( x(t),\dot{x}(t) ,t \right) \right) & =  e^{\eta t}  \sum_{i=1}^{n}{ g_{ik}(x(t)) \frac{dv^i}{dt}(t)}   +e^{\eta t}  \sum_{i,j=1}^{n}{ \frac{\partial g_{kj}}{\partial q^i}(x(t)) v^i(t) v^j(t)  } \\ & \qquad \qquad \quad\qquad + \eta e^{\eta t} \sum_{i=1}^{n}{g_{ik}(x(t)) v^i(t)}, 
	\end{align*}
	\[  \frac{\partial \mathcal{L}^{SC}}{\partial q^k}  \left( x(t),\dot{x}(t) ,t \right)  =  e^{\eta t} \sum_{i,j=1}^{n}{ \frac{\partial g_{ij}}{\partial q^k}(x(t)) v^i(t)v^j(t)} - e^{\eta t} \frac{\partial f}{\partial q^k} (x(t)). \]
	If we multiply both terms by $e^{-\eta t}$, the Euler--Lagrange equations~\eqref{eq: EL Basic} for the Lagrangian $\mathcal{L}^{SC}$ are given, for $k=1, \ldots , n$, by
	\begin{align*}  0 = \sum_{i=1}^{n}{ g_{ik}(x(t)) \frac{dv^i}{dt}(t)}  & +   \sum_{i,j=1}^{n}{ \frac{\partial g_{kj}}{\partial q^i}(x(t)) v^i(t) v^j(t)  } +  \eta  \sum_{i=1}^{n}{g_{ik}(x(t)) v^i(t)} \\ &  \qquad\quad -  \frac{1}{2} \sum_{i,j=1}^{n}{ \frac{\partial g_{ij}}{\partial q^k}(x(t)) v^i(t)v^j(t)} +  \frac{\partial f}{\partial q^k} (x(t)) .
	\end{align*}
	Rearranging terms, and multiplying by the matrix $(g^{ij})$ which is the inverse of $(g_{ij})$, we get, for $k=1,\ldots n$,  the equation
	\[ \left( \frac{dv^k}{dt}(t)   +   \sum_{i,j=1}^{n}{ \Gamma^{k}_{ij}(x(t)) v^i(t) v^j(t)  }  \right)  + \eta v^k(t) +  \left( \text{gradf}(x(t)) \right)^k  = 0,\]
	where $\Gamma^{k}_{ij}$ are the Christoffel symbols given by
	$ \Gamma^{k}_{ij} = \frac{1}{2} \sum_{l=1}^{n}{ g^{kl} \left[ \frac{\partial g_{jl}}{\partial x^i} + \frac{\partial g_{li}}{\partial x^j} - \frac{\partial g_{ij}}{\partial x^l}\right]},$ which gives the desired Euler--Lagrange equation. \qed }
\end{theorem}

\section{Proof of the Convergence Rates}\label{Section: Appendix Proof Convergence Rate}

The proofs of the convergence rates of solutions to the Bregman Euler--Lagrange equations are inspired by those of Theorems 5 and 6 from \cite{alimisis2020}, and make use of Lemmas 2 and 12 therein:

\begin{lemma} \label{Lemma 2}
	Given a Riemannian manifold $\mathcal{Q}$ with sectional curvature bounded above by $K_{\max}$ and below by $K_{\min}$, with $\zeta$ given by equation~\eqref{eq: zeta}, and such that
	\begin{align*} 
		\emph{diam}(\mathcal{Q}) <
		\begin{cases}
			\frac{\pi}{\sqrt{K_{\max}}}  & \quad \emph{if }  K_{\max} > 0 \\ \infty   & \quad \emph{if }  K_{\max} \leq 0 
		\end{cases}  	,
	\end{align*} 
	we have that
	\[ \langle  \nabla_{\dot{X}} \emph{Log}_X(p) , -\dot{X} \rangle \leq \zeta \| \dot{X}\|^2. \] 
\end{lemma}  
\begin{lemma} \label{Lemma 12}
	Given a point $q$ and a smooth curve $X(t)$ on a Riemannian manifold $\mathcal{Q}$,
	\begin{align*} \frac{d}{dt} \| \emph{Log}_{X(t)} (q) \|^2 = 2\langle \emph{Log}_{X(t)}(q) , \nabla_{\dot{X}} \emph{Log}_{X(t)}(q)  \rangle = 2 \langle \emph{Log}_{X(t)}(q) , -\dot{X}(t) \rangle. \\
	\end{align*} 
\end{lemma}

\begin{theorem}
	Suppose $f : \mathcal{Q} \rightarrow \mathbb{R}$ is a $\lambda$-weakly-quasi-convex function, and suppose that Assumption~\ref{Assumption 1} is satisfied. Then, any solution $X(t)$ of the Bregman Euler--Lagrange equation
	\[
	\nabla_{\dot{X}}\dot{X}  +\left( \lambda^{-1} \zeta e^{\alpha_t} -  \dot{\alpha}_t \right) \dot{X} + e^{2\alpha_t+\beta_t }  \emph{gradf}(X) = 0,
	\] 
	with $X(0)=x_0$ and $\dot{X}(0) = 0$, converges to a minimizer $x^*$ of $f$ with rate
	\[
	f(X(t)) - f(x^*) \leq \frac{2 \lambda^2 e^{\beta_0} \left( f(x_0) - f(x^*) \right) +  \zeta \| \emph{Log}_{x_0}{(x^*)} \|^2 }{2\lambda^2 e^{\beta_t}} .
	\] 
	\proof{ Let
		\[ \mathcal{E}(t) = \lambda^2 e^{\beta_t} \left( f(X) - f(x^*) \right)  + \frac{1}{2} (\zeta -1)\| \text{Log}_{X}(x^*) \|^2	 + \frac{1}{2} \bigg\| \lambda e^{-\alpha_t} \dot{X} -\text{Log}_{X}(x^*)  \bigg\|^2. \]
		Then, using Lemma~\ref{Lemma 12},
		\begin{align*}
			\dot{\mathcal{E}}(t)  & =    \lambda^2 \dot{\beta}_t e^{\beta_t} \left( f(X) - f(x^*) \right) + \lambda^2 e^{\beta_t} \langle \text{gradf}(X), \dot{X}\rangle  + (\zeta -1) \langle \text{Log}_{X}(x^*) , -\dot{X} \rangle \\ &  \qquad \qquad \qquad  \qquad +   \langle \lambda e^{-\alpha_t} \dot{X} -\text{Log}_{X}(x^*)   ,  -\dot{\alpha}_t \lambda e^{-\alpha}  \dot{X} + \lambda e^{-\alpha_t} \nabla_{\dot{X}} \dot{X}    -\nabla_{\dot{X}}  \text{Log}_{X}(x^*)   \rangle  \\ & = \lambda^2 \dot{\beta}_t e^{\beta_t} \left( f(X) - f(x^*) \right) + \lambda^2 e^{\beta_t} \langle \text{gradf}(X), \dot{X}\rangle + (\zeta -1) \langle \text{Log}_{X}(x^*) , -\dot{X} \rangle \\ & \qquad \qquad  \qquad \qquad +  \langle \lambda e^{-\alpha_t} \dot{X} -\text{Log}_{X}(x^*)   , \lambda e^{-\alpha_t} \left(-\dot{\alpha}_t  \dot{X} + \nabla_{\dot{X}} \dot{X}  \right)  -\nabla_{\dot{X}}  \text{Log}_{X}(x^*)   \rangle .
		\end{align*} 
		Now, from the Bregman Euler--Lagrange equation,
		\[        -\dot{\alpha}_t \dot{X} +  \nabla_{\dot{X}} \dot{X}   = - \lambda^{-1} \zeta e^{\alpha_t} \dot{X}  -  e^{2\alpha_t+\beta_t }  \text{gradf}(X).         \]  
		Thus,
		\begin{align*}
			\dot{\mathcal{E}}(t)   & = \lambda^2 \dot{\beta_t} e^{\beta_t} \left( f(X) - f(x^*) \right) + \lambda^2 e^{\beta_t} \langle \text{gradf}(X), \dot{X}\rangle + (\zeta -1) \langle \text{Log}_{X}(x^*) , -\dot{X} \rangle \\ & \qquad \qquad  \qquad  \qquad  +  \langle \lambda e^{-\alpha_t} \dot{X} -\text{Log}_{X}(x^*)   ,  - \zeta \dot{X}  - \lambda  e^{\alpha_t+\beta_t}  \text{gradf}(X)     -\nabla_{\dot{X}}  \text{Log}_{X}(x^*)   \rangle  \\  & =  \lambda^2 \dot{\beta}_t e^{\beta_t} \left( f(X) - f(x^*) \right) + \lambda^2 e^{\beta_t} \langle \text{gradf}(X), \dot{X}\rangle  + (\zeta -1) \langle \text{Log}_{X}(x^*) , -\dot{X} \rangle   - \lambda \zeta e^{-\alpha_t} \langle  \dot{X}   ,  \dot{X}  \rangle  \\ & \qquad  \qquad \qquad  \qquad  -\lambda^2  e^{\beta_t}  \langle \dot{X}   ,     \text{gradf}(X)     \rangle - \lambda  e^{-\alpha_t} \langle  \dot{X}   ,   \nabla_{\dot{X}}  \text{Log}_{X}(x^*)   \rangle  +  \zeta \langle  \text{Log}_{X}(x^*)   ,   \dot{X}    \rangle  \\& \qquad \qquad \qquad \qquad + \lambda  e^{\alpha_t+\beta_t }  \langle  \text{Log}_{X}(x^*)   ,   \text{gradf}(X)      \rangle +  \langle  \text{Log}_{X}(x^*)   ,  \nabla_{\dot{X}}  \text{Log}_{X}(x^*)   \rangle.
		\end{align*} 
		Canceling the $ \langle \text{gradf}(X), \dot{X}\rangle $ and $\langle \text{Log}_{X}(x^*) , -\dot{X} \rangle$ terms out using Lemma~\ref{Lemma 12}, we get
		\begin{align*}
			\dot{\mathcal{E}}(t)  & = \lambda^2 \dot{\beta}_t e^{\beta_t} \left( f(X) - f(x^*) \right)    +  \lambda e^{\alpha_t+\beta_t }  \langle  \text{Log}_{X}(x^*)   ,   \text{gradf}(X)      \rangle  \\
			& \qquad \qquad \qquad \qquad \qquad \qquad    - \lambda \zeta e^{-\alpha_t} \langle  \dot{X}   ,  \dot{X}  \rangle  - \lambda  e^{-\alpha_t} \langle  \dot{X}   ,   \nabla_{\dot{X}}  \text{Log}_{X}(x^*)   \rangle  \\ & = \lambda e^{\beta_t} \left[ \dot{\beta}_t  \lambda \left( f(X) - f(x^*) \right)  + e^{\alpha_t} \langle  \text{Log}_{X}(x^*)   ,   \text{gradf}(X)      \rangle  \right]  \\ & \qquad \qquad \qquad \qquad \qquad \qquad  - \lambda  e^{-\alpha_t} \left[ \zeta\langle  \dot{X}   ,  \dot{X}  \rangle  +  \langle  \dot{X}   ,   \nabla_{\dot{X}}  \text{Log}_{X}(x^*)   \rangle  \right].
		\end{align*} 
		Now, since $f$ is geodesically $\lambda$-weakly-quasi-convex, we have that  
		\[  \lambda \left( f(X) - f(x^*) \right) + \langle \text{Log}_{X}(x^*) , \text{gradf}(X) \rangle  \leq 0, \]
		so the ideal scaling equation $\dot{\beta}_t \leq e^{\alpha_t}$ implies that 
		\[   \lambda  e^{\beta_t} \left[ \dot{\beta}_t \lambda \left( f(X) - f(x^*) \right)  + e^{\alpha_t} \langle  \text{Log}_{X}(x^*)   ,   \text{gradf}(X)      \rangle  \right]     \leq 0. \] 
		Moreover, Lemma~\ref{Lemma 2} yields $   \left[ \zeta \langle \dot{X} , \dot{X} \rangle  + \langle \dot{X}  ,  \nabla_{\dot{X}}  \text{Log}_{X}(x^*)   \rangle \right] \geq 0$, so \[  - \lambda e^{-\alpha_t} \left[ \zeta\langle  \dot{X}   ,  \dot{X}  \rangle  +  \langle  \dot{X}   ,   \nabla_{\dot{X}}  \text{Log}_{X}(x^*)   \rangle  \right] \leq 0. \] Therefore, $\dot{\mathcal{E}}(t) \leq 0$, and so
		\begin{align*}
			\lambda^2 e^{\beta_t} \left( f(X) - f(x^*) \right)  & \leq \lambda^2 e^{\beta_t} \left( f(X) - f(x^*) \right)  + \frac{1}{2} (\zeta -1)\| \text{Log}_{X}(x^*) \|^2	 + \frac{1}{2} \bigg\| \lambda e^{-\alpha_t} \dot{X} -\text{Log}_{X}(x^*)  \bigg\|^2  \\ & = \mathcal{E}(t)  \leq \mathcal{E}(0) =  \lambda^2 e^{\beta_0} \left( f(x_0) - f(x^*) \right)  + \frac{1}{2} \zeta \| \text{Log}_{x_0}(x^*) \|^2,
		\end{align*}
		which gives the desired rate of convergence
		\[    
		f(X(t)) - f(x^*) \leq \frac{2 \lambda^2 e^{\beta_0} \left( f(x_0) - f(x^*) \right) +  \zeta \| \text{Log}_{x_0}{(x^*)} \|^2 }{2 \lambda^2 e^{\beta_t}}.
		\] \qed }
\end{theorem}

\section{Proof of Existence Theorems} 

\subsection{Convex and Weakly-Quasi-Convex Cases}\label{Section: Appendix Existence}

\begin{theorem} 
	Suppose Assumption~\ref{Assumption 1} is satisfied, and let $C, p >0$ and $v>1$ be given constants. Then the differential equation
	\[
		\nabla_{\dot{X}}\dot{X}  + \frac{v}{t} \dot{X} + Ct^{p-2} \emph{gradf}(X) = 0,
	\] 
	has a global solution $X:[0,\infty)\rightarrow \mathcal{Q}$ under the initial conditions $X(0)=x_0 \in \mathcal{Q}$ and $\dot{X}(0)=0.$
 \proof{The proof is similar to that of Lemma 3 in \cite{alimisis2020}, which extended Theorem 1 in \cite{SuBoCa16} to the Riemannian setting. We first define a family of smoothed equations for which we then show existence of a solution for all time. After choosing an equicontinuous and uniformly bounded subfamily of smoothed solutions, we use the Arzela--Ascoli Theorem on the complete Riemannian manifold $\mathcal{Q}$ to obtain a subsequence converging uniformly, and argue that the limit of this subsequence solves the original problem. When $p=2$, we recover the simpler case considered in Lemma 3 of \cite{alimisis2020}, so we assume $p \neq 2$ in this proof. Consider the following families of smoothed equations for $\delta>0$:
		\begin{align*}
			&	\nabla_{\dot{X}}\dot{X}  + \frac{v}{\max{(\delta,t)}} \dot{X} + C (\max{(\delta,t)})^{p-2} \text{gradf}(X) = 0    &  \text{if } p<2, \\
			&	\nabla_{\dot{X}}\dot{X}  + \frac{v }{\max{(\delta,t)}} \dot{X} + C t^{p-2} \text{gradf}(X) = 0     & \text{if } p>2.
		\end{align*}
		Exp and Log are defined globally on $\mathcal{Q}$ by Assumption~\ref{Assumption 1}, so we can choose geodesically normal coordinates $\phi = \psi^{-1}$ around $x_0$ defined globally on $\mathcal{Q}$ and put $c= \phi \circ X$. Using the smoothness of $f$ and letting $u=\dot{c}$ gives a system of first-order ODEs defining a local representation for a vector field in $T\mathcal{Q}$, and Section IV.3 of \cite{Lang1999} guarantees that the smoothed ODE has a unique solution $X_\delta$ locally around 0. Actually, $X_\delta$ exists on $[0,\infty)$. Indeed, by contradiction, let $[0,T)$ be the maximal interval of existence of $X_\delta$, for some finite $T>0$. Using \[\frac{d}{dt}f(X_\delta(t))  = \langle \text{gradf}(X_\delta), \dot{X}_\delta \rangle \] gives 
			\small	\begin{align*}
			\frac{d}{dt}f(X_\delta) & =  - \frac{\delta^{2-p}}{C} \langle \nabla_{\dot{X}_\delta}\dot{X}_\delta, \dot{X}_\delta \rangle - \frac{v \delta ^{1-p}}{C} \langle \dot{X}_\delta, \dot{X}_\delta \rangle = - \frac{\delta^{2-p}}{2C} \frac{d}{dt} \| \dot{X}_\delta \|^2 - \frac{v \delta ^{1-p}}{C} \| \dot{X}_\delta \|^2 \quad  \text{ }\text{if } \delta>t, \text{ } p<2, & \\ \frac{d}{dt}f(X_\delta)  & =  - \frac{t^{2-p}}{C} \langle \nabla_{\dot{X}_\delta}\dot{X}_\delta, \dot{X}_\delta \rangle - \frac{v t ^{2-p}}{C\delta } \langle \dot{X}_\delta, \dot{X}_\delta \rangle = - \frac{t^{2-p}}{2C} \frac{d}{dt} \| \dot{X}_\delta \|^2 - \frac{v t ^{2-p}}{C\delta } \| \dot{X}_\delta \|^2 \qquad \text{if } \delta>t, \text{ } p>2 ,&\\
			\frac{d}{dt}f(X_\delta) & =  - \frac{t^{2-p}}{C} \langle \nabla_{\dot{X}_\delta}\dot{X}_\delta, \dot{X}_\delta \rangle - \frac{v t^{1-p}}{C} \langle \dot{X}_\delta, \dot{X}_\delta \rangle   = -\frac{1}{2C} \frac{d}{dt} \left( t^{2-p} \|  \dot{X}_\delta \|^2 \right) - \frac{2v(2-p)-1}{2C(2-p)} t^{1-p} \|  \dot{X}_\delta \|^2  \quad \text{if } \delta<t.& 
		\end{align*}
	\normalsize	Let $\theta = \frac{2v(2-p)-1}{2C(2-p)} $. Integrating and using the Cauchy-Schwarz inequality for the $p<2$ case gives
\small	\begin{align*}
	& \int_{0}^{T}{  \sqrt{(\max{(\delta , t)})^{1-p}} \| \dot{X}_\delta  \| dt} = \int_{0}^{\delta}{  \sqrt{\delta^{1-p}} \| \dot{X}_\delta  \| dt}  + \int_{\delta}^{T}{  \sqrt{t^{1-p}} \| \dot{X}_\delta  \| dt} 	 \\ & \quad \qquad \leq   \sqrt{  \frac{C\delta }{v} (f(x_0)-\inf_{u}{f(u)}) + \frac{\delta^{2-p} }{2 v} \left(  \| \dot{X}_\delta (0)\|^2 -  \inf_{t\in [0,T) }{\| \dot{X}_\delta (t)\|^2}  \right)}   \\ &   \qquad \qquad \qquad + \sqrt{  \frac{T-\delta }{\theta} (f(X_\delta (\delta ))-\inf_{u}{f(u)}) + \frac{T-\delta  }{2C\theta} \left( \delta^{2-p} \| \dot{X}_\delta (\delta )\|^2 -  \inf_{t\in [0,T) }{t^{2-p} \| \dot{X}_\delta (t)\|^2}  \right)}    < \infty,
\end{align*}
	\normalsize	since $f$ is bounded below by Assumption~\ref{Assumption 1}. If $\delta \geq T$, then $ \sqrt{\delta^{1-p}} \dot{X}_\delta $ is integrable on $[0,T)$. If $\delta< T$, then the integrals on $[0,T)$ and $[0,\delta)$ are finite, so the integral on $[\delta ,T)$ must also be finite, and thus $ \sqrt{t^{1-p}}  \dot{X}_\delta $ is integrable on $[\delta , T)$. Now, $\| \int_{a}^{T}{\dot{X}_\delta dt}\| \leq \int_{a}^{T}{\| \dot{X}_\delta \| dt} <\infty$ for $a=0,\delta$ implies that $\lim_{t\rightarrow T}{X_\delta(t)}$ exists. Since $\mathcal{Q}$ is complete by Assumption~\ref{Assumption 1}, the limit is in $\mathcal{Q}$, contradicting the maximality of $[0,T)$. The $p>2$ case is similar: the integrand is replaced by $\sqrt{t^{2-p}(\max{(\delta , t)})^{-1}} \| \dot{X}_\delta  \|$, and the integral on $[\delta,T)$ remains unchanged while the integral on $[0,\delta)$ can be bounded by the same expression using $t<\delta$. Thus, in both cases, we can find a solution $X_\delta : [0,\infty ) \rightarrow \mathcal{Q}$ to the smooth initial-value ODE, and its corresponding solution  $X_\delta : [0,\infty ) \rightarrow \mathbb{R}^n$ in local coordinates. 
	
	Now define \[M_\delta (t) = \sup_{u\in (0,t] }{\frac{\| \dot{X}_\delta (u) \|}{u} }. \] 
When $0 < t \leq \delta,$ the smoothed ODE can be written as
		\begin{align*}
			\nabla_{\dot{X}_\delta} \left( 	\dot{X} _\delta e^{\frac{v}{\delta} }  \right)  = -  C\delta^{p-2} \text{gradf}(X_\delta) e^{\frac{v}{\delta} }   \text{  if } p< 2,
			\quad 	\nabla_{\dot{X}_\delta} \left( 	\dot{X} _\delta e^{\frac{v}{\delta} }  \right)  = -  C t^{p-2} \text{gradf}(X_\delta) e^{\frac{v}{\delta} }   \text{  if } p> 2.
		\end{align*}
		Thus, we can use Lemma 4 in \cite{alimisis2020} to get for $p > 2$ that 
		\begin{align*}
			\Gamma_{X_\delta (t)}^{x_0} \dot{X}_\delta (t) & = -e^{-\frac{v}{\delta} t }    \int_{0}^{t}{  \left( \Gamma_{X_\delta (u)}^{x_0}  \text{gradf}(X_\delta (u)) -   \Gamma_{X_\delta (u)}^{x_0}  \Gamma(X_\delta)_{x_0}^{X_\delta (u)}  \text{gradf}(x_0) \right) Cu^{p-2} e^{\frac{v}{\delta}u } du }   \\ & \qquad  \qquad\qquad      -e^{-\frac{v}{\delta} t }  \int_{0}^{t}{   Cu^{p-2}  \Gamma_{X_\delta (u)}^{x_0}  \Gamma(X_\delta)_{x_0}^{X_\delta (u)}  \text{gradf}(x_0)    e^{\frac{v}{\delta}u } du  }.
		\end{align*} 
From the Lipschitz assumption on $f$, we have that
		\[ \big\| \text{gradf}(X_\delta (u)) - \Gamma_{x_0}^{X_\delta (u)} \text{gradf}(x_0) \big\|  \leq L \int_{0}^{u}{\| \dot{X}_{\delta}(s) \| ds} = L \int_{0}^{u}{s \frac{ \| \dot{X}_{\delta}(s)\|}{s}  ds}  \leq \frac{1}{2} LM_\delta (u) u^2.\]
		Thus, since parallel transport preserves inner products,
		\begin{align*}
			\small \frac{\|\dot{X}_\delta(t) \|}{t}   & \leq   \left(    \frac{1}{2} C L M_{\delta}(\delta) \delta^{p} + C\delta^{p}  \| \text{gradf}(x_0)\|  \right) \frac{e^{-\frac{v}{\delta} t } }{t}   \int_{0}^{t}{  e^{\frac{v}{\delta}u } du  }   \\ & \leq    \left(    \frac{1}{2} C L M_{\delta}(\delta) \delta^{p} + C\delta^{p}  \| \text{gradf}(x_0)\|  \right)  \frac{\delta}{vt}(1-e^{-\frac{v}{\delta} t }) \leq  \frac{1}{2} C L M_{\delta}(\delta) \delta^{p} + C\delta^{p}  \| \text{gradf}(x_0)\|. 
		\end{align*} 
		Taking the supremum over $0<t \leq \delta$ and rearranging gives for $\delta < \delta_M = \left( \frac{2}{CL} \right)^{\frac{1}{p}}$ that \[ M_\delta (\delta) \leq  \frac{ 2C \delta^{p} \| \text{gradf}(x_0) \|}{2-CL \delta^{p}}.\]
		The case $p < 2$ is done exactly in the same way except that we do not need to bound $u^{p-2}$ by $\delta^{p-2}$ in the integrals since the $t^{p-2}$ term in the differential equation is already replaced by $\delta^{p-2}$. 
		
		\noindent Note that when $\delta < \delta_M$ and $\delta < t  < t_M = \left( \frac{2(v+p+1)}{CL} \right)^{\frac{1}{p}}$,  the smoothed ODE can be rewritten as
		\[\frac{d}{dt} \left( 	t^v \dot{X} _\delta(t)   \right)  = -  Ct^{v+p-2} \text{gradf}(X_\delta).\]
		Therefore, we can use Lemma 4 in \cite{alimisis2020} once again to obtain
		\begin{align*}
			\Gamma_{X_\delta (t)}^{X_\delta(\delta)} t^v \dot{X}_\delta (t) - \delta^v \dot{X}_\delta (\delta) & =  \int_{0}^{t}{  \left( \Gamma_{X_\delta (u)}^{X_\delta(\delta)}  \text{gradf}(X_\delta (u)) -   \Gamma_{X_\delta (u)}^{X_\delta(\delta)}  \Gamma(X_\delta)_{x_0}^{X_\delta (u)}  \text{gradf}(x_0) \right) Cu^{v+ p-2}  du }   \\ & \qquad  \qquad\qquad      -  \int_{0}^{t}{   Cu^{v+p-2}  \Gamma_{X_\delta (u)}^{X_\delta(\delta)}  \Gamma(X_\delta)_{x_0}^{X_\delta (u)}  \text{gradf}(x_0)  du  }.
		\end{align*} 
		Using the fact that parallel transport preserves inner products, and dividing by $t^{v+1}$ gives
		\begin{align*}
			\small \frac{\|\dot{X}_\delta(t) \|}{t}   & \leq  \frac{\delta^{v+1}}{t^{v+1}} \frac{\|\dot{X}_\delta(\delta) \|}{\delta} + \frac{CL}{2t^{v+1}} \int_{\delta}^{t}{M_\delta(u) u^{v+p} du} + \frac{C}{t^{v+1}} \| \text{gradf}(x_0)\| \int_{\delta}^{t}{u^{v+p-2} du}
			\\ & \leq    \frac{\delta^{v+1}}{t^{v+1}} \frac{ 2C \delta^{p} \| \text{gradf}(x_0) \|}{2-CL \delta^{p}} + \frac{CL}{2(v+p+1)} M_\delta (t) t^p + \frac{C(t^{v+p-1} - \delta^{v+p-1})}{(v+p-1)t^{v+1}} \| \text{gradf}(x_0)\| ,
		\end{align*} 
		and since this upper bound is an increasing function of $t$, we have for any $t'\in (\delta,t) $ that 
		\small \begin{align*}
			\small \frac{\|\dot{X}_\delta(t') \|}{t'}   &  \leq    \frac{ 2C \delta^{p} \| \text{gradf}(x_0) \|}{2-CL \delta^{p}} + \frac{CL}{2(v+p+1)} M_\delta (t) t^p + \frac{Ct^{p-2} }{v+p-1} \| \text{gradf}(x_0)\|.
		\end{align*} 
	\normalsize	Taking the supremum over all $t' \in (0,t)$ gives for $\delta < \delta_M$ and $\delta < t <t_M$, 
		 \begin{align*}
			M_\delta (t) \leq    \frac{1}{1-\frac{CL}{2(v+p+1)} t^p} \left( \frac{ 2C \delta^{p} }{2-CL \delta^{p}} + \frac{Ct^{p-2} }{v+p-1}  \right)   \| \text{gradf}(x_0) \|.
		\end{align*} 
		Now consider the family of functions
		\[
			\mathcal{F} = \Big\{  X_\delta :\left[ 0 , T  \right] \rightarrow \mathbb{R} \big| \delta = 2^{-n}\tilde{\delta} , n=0,1,\dots    \Big\} ,
		\]
		where $T = \left( \frac{v+p+1}{CL} \right)^{\frac{1}{p}}$ and $\tilde{\delta} =\left( \frac{1}{CL} \right)^{\frac{1}{p}} $.
		By definition of $M_\delta$, we have for $t\in [ 0 , T]  $ and $\delta \in (0,  \tilde{\delta})$  that
		\small \[ \| \dot{X}_\delta \| \leq  T  M_\delta (T) \leq    2CT \left( \tilde{\delta}+  \frac{CT^{p-2}}{v+p-1}\right) \quad \text{and} \quad  d(X_\delta (t) ,X_\delta (0)) \leq \int_{0}^{t}{\| \dot{X}_\delta (u) \| du} \leq t  \| \dot{X}_\delta \|  \leq  T\| \dot{X}_\delta \| . \]  
		\normalsize Thus, $\mathcal{F}$ is equicontinuous and uniformly bounded, and the Riemannian manifold $\mathcal{Q}$ is complete by Assumption~\ref{Assumption 1}, so by the Arzela--Ascoli Theorem (Theorem 17 in \cite{Kelley1975}), $\mathcal{F}$ contains a subsequence that converges uniformly on $[0,T]$ to some function $X^*$. The same argument as in part 5 of the proof of Lemma 3 of \cite{alimisis2020} shows that $X^*$ is a solution to the original initial-value ODE on $[0,T]$ which can then be extended to get a global solution on $[0,\infty)$. \qed 
	}
\end{theorem}

\subsection{Strongly Convex Case}\label{Section: Appendix Existence SC}

\begin{theorem} 
Suppose that Assumption~\ref{Assumption 1} is satisfied, and that $\eta  >0$ is a given constant. Then, the differential equation
\[
	\nabla_{\dot{X}}\dot{X}  + \eta  \dot{X} + \emph{gradf}(X) = 0,
\] 
has a global solution $X:[0,\infty)\rightarrow \mathcal{Q}$ under the initial conditions $X(0)=x_0 \in \mathcal{Q}$ and $\dot{X}(0)=0.$
\proof{Exp and Log are defined globally on $\mathcal{Q}$ by Assumption~\ref{Assumption 1}, so we can choose geodesically normal coordinates $\phi = \psi^{-1}$ around $x_0$ defined globally on $\mathcal{Q}$ and put $c= \phi \circ X$. As in \cite{alimisis2020}, using the smoothness of $f$ and letting $u=\dot{c}$ gives a system of first-order ODEs which defines a local representation for a vector field in $T\mathcal{Q}$, and results from Section IV.3 of \cite{Lang1999} guarantee that the initial-value differential equation has a unique solution locally around 0. It remains to show that this solution actually exists on $[0,\infty)$. Towards contradiction, suppose $[0,T)$ is the maximal interval of existence of the solution $X$, for some finite $T>0$. Then,  
	\[\frac{d}{dt}f(X(t)) = \langle \text{gradf}(X), \dot{X} \rangle = - \langle \nabla_{\dot{X}}\dot{X}, \dot{X} \rangle - C \langle \dot{X}, \dot{X} \rangle = -\frac{1}{2} \frac{d}{dt} \| \dot{X} \|^2  - C \| \dot{X} \|^2.\] 
	Rearranging, integrating both sides and using the Cauchy-Schwarz inequality gives
	\[\int_{0}^{T}{ \| \dot{X} \| dt} = \sqrt{T(f(x_0)-\inf_{u}{f(u)}) + \frac{T}{2} \left( \| \dot{X}(0)\|^2 - \inf_{t\in [0,T) }{ \| \dot{X}(t)\|^2}  \right)} < \infty,\]
	since $f$ is bounded from below by Assumption~\ref{Assumption 1}. Thus, $\lim_{t\rightarrow T}{X(t)}$ exists, and since $\mathcal{Q}$ is complete, the limit is in $\mathcal{Q}$, contradicting the maximality of $[0,T)$. This completes the proof.  \qed }
\end{theorem} \hfill

\section{Proof of Invariance Theorem}  \label{Section: Appendix Proof Invariance}

\begin{theorem} 
	Suppose that Assumption~\ref{Assumption 1} is satisfied and that the curve $X(t)$ satisfies the Riemannian Bregman Euler--Lagrange equation~\eqref{eq: EL} corresponding to $\mathcal{L}_{\alpha,\beta,\gamma}$. Then the reparametrized curve $ X(\tau(t))$ satisfies the Bregman Euler--Lagrange equation~\eqref{eq: EL} corresponding to the modified Riemannian Bregman Lagrangian $\mathcal{L}_{\tilde{\alpha},\tilde{\beta},\tilde{\gamma}}$ where $\tilde{\alpha}_t = \alpha_{\tau(t)} + \log{\dot{\tau}(t) }$,  $\tilde{\beta}_t = \beta_{\tau(t)}$, and $\tilde{\gamma}_t = \gamma_{\tau(t)}$. Furthermore $\alpha,\beta,\gamma$ satisfy the ideal scaling conditions (\ref{IdealScaling}) if and only if $\tilde{\alpha},\tilde{\beta},\tilde{\gamma}$ do.
	\proof{ Let $Y(t) = X(\tau (t)) $. Then
		\[ \dot{Y}(t) = \dot{\tau}(t) \dot{X}(\tau(t)), \qquad \text{and} \qquad \nabla_{\dot{Y}(t)}\dot{Y}(t)  = \ddot{\tau}(t) \dot{X}(\tau(t)) + \dot{\tau}^2(t) \nabla_{\dot{X}(\tau(t))}\dot{X}(\tau(t)).\]
		Inverting these relations gives
		\[\dot{X}(\tau(t)) =  \frac{1}{ \dot{\tau}(t) }\dot{Y}(t),  \qquad \text{and} \qquad \nabla_{\dot{X}(\tau(t))}\dot{X}(\tau(t)) = \frac{1}{ \dot{\tau}^2(t) } \nabla_{\dot{Y}(t)}\dot{Y}(t)  - \frac{\ddot{\tau}(t)}{\dot{\tau}^3(t)} \dot{Y}(t).\]
		The Bregman Euler--Lagrange equation~\eqref{eq: EL} at time $\tau(t)$ is given by
		\[   \nabla_{\dot{X}(\tau(t)) }\dot{X}(\tau(t))  +\left( \lambda^{-1} \zeta e^{\alpha_{\tau(t)}} -  \dot{\alpha}_{\tau(t)} \right) \dot{X}(\tau(t))+ e^{2\alpha_{\tau(t)}+\beta_{\tau(t)} }  \text{gradf}(X(\tau(t))) = 0.   \]
		Substituting the expressions for $X(\tau(t)), \dot{X}(\tau(t))$ and $\nabla_{\dot{X}(\tau(t))}\dot{X}(\tau(t))$ in terms of $Y(t)$ and its derivatives, and multiplying by $\dot{\tau}^2(t)$, we get
		\[    \nabla_{\dot{Y}(t)}\dot{Y}(t)  - \frac{\ddot{\tau}(t)}{\dot{\tau}(t)} \dot{Y}(t) +\left( \lambda^{-1} \zeta e^{\alpha_{\tau(t)}} -  \dot{\alpha}_{\tau(t)} \right)  \dot{\tau}(t) \dot{Y}(t) + \dot{\tau}^2(t) e^{2\alpha_{\tau(t)}+\beta_{\tau(t)} }  \text{gradf}(Y(t)) = 0.   \]
		Substituting the expressions for $\alpha, \beta, \gamma$ in terms of $\tilde{\alpha}, \tilde{\beta}, \tilde{\gamma}$ yields
		\[    \nabla_{\dot{Y}(t)}\dot{Y}(t)  - \frac{\ddot{\tau}(t)}{\dot{\tau}(t)} \dot{Y}(t) +\left( \lambda^{-1} \zeta \frac{1}{\dot{\tau}(t)} e^{\tilde{\alpha}_{t}} -  \frac{1}{\dot{\tau}(t)}  \left[ \dot{\tilde{\alpha}} (t) + \frac{\ddot{\tau}(t)}{\dot{\tau}(t)} \right]  \right)  \dot{\tau}(t) \dot{Y}(t) + e^{2\tilde{\alpha}_{t}+\tilde{\beta}_{t} }  \text{gradf}(Y(t)) = 0.   \] 
		This gives the Bregman Euler--Lagrange equation~\eqref{eq: EL} corresponding to $\mathcal{L}_{\tilde{\alpha},\tilde{\beta},\tilde{\gamma}}$,
		\[    \nabla_{\dot{Y}(t)}\dot{Y}(t)   +\left( \lambda^{-1} \zeta  e^{\tilde{\alpha}_{t}} -  \frac{1}{\dot{\tau}(t)}  \dot{\tilde{\alpha}} (t)  \right) \dot{Y}(t) + e^{2\tilde{\alpha}_{t}+\tilde{\beta}_{t} }  \text{gradf}(Y(t)) = 0.   \] 
		The fact that the parameters $\alpha,\beta,\gamma$ satisfy the ideal scaling conditions (\ref{IdealScaling}) if and only if the parameters $\tilde{\alpha},\tilde{\beta},\tilde{\gamma}$ do is established in the proof of Theorem 1.2 of \cite{WiWiJo16}.	\qed} \\
\end{theorem}

\bibliography{RiemannVariational}
\bibliographystyle{siamplain}

\end{document}